\documentclass[12pt]{article}

\usepackage{amsmath,amsthm,amsfonts}

\parskip=\medskipamount
\def\cinfty#1{C^{\scriptscriptstyle\infty}(#1)}
\def\vectorfields#1{{\cal X}(#1)}
\def\ov#1{\overline{#1}}

\def\fpd#1#2{\frac{\partial #1}{\partial #2}}
\def\R{{\rm I\kern-.20em R}}
\def\sode{{\sc Sode}}
\def\proof{{\sc Proof.}}
\newtheorem{prop}{\bf Proposition}
\newtheorem{lemma}{\bf Lemma}

\newtheorem{dfn}{\bf Definition}
\newtheorem{cor}{\bf Corollary}

\def\base#1{{\ov e}_{#1}}
\def\cosyst{(e_0; \{ \base \alpha \})}
\def\baseX#1{{\mathcal X}_{#1}}
\def\baseV#1{{\mathcal V}_{#1}}
\def\baseH#1{{\mathcal H}_{#1}}

\def\cans{{\mathcal I}}

\def\H#1{{#1}^{\scriptscriptstyle H}}
\def\V#1{{#1}^{\scriptscriptstyle V}}
\def\subH#1{{#1}_{\scriptscriptstyle H}}
\def\subV#1{{#1}_{\scriptscriptstyle V}}

\def\tpi{\tilde\pi}
\def\trho{\tilde\rho}

\def\prol#1#2{{T^{#1}{#2}}}
\def\prolr#1{\prol{\varrho}{#1}}
\def\proltr#1{\prol{\tilde{\rho}}{#1}}

\newcommand{\mybox}[1]{\makebox(0,0){\footnotesize{#1}}}
\newlength{\savelen}

\def\cM{c_{\scriptscriptstyle M}}
\def\psiM{\psi_{\scriptscriptstyle M}}

\def\vl{v_{\cM^*E}}

\def\cd#1{{\rm D}_{#1}}
\def\cdtilde#1{{\tilde{\rm D}}_{#1}}
\def\cdhat#1{{\hat{\rm D}}_{#1}}
\def\cdov#1{{\ov{\rm D}}_{#1}}

\def\cdhattilde#1{{\tilde{{\hat{\rm D}}}}_{#1}}
\def\cdhatov#1{{\ov{{\hat{\rm D}}}}_{#1}}

\begin{document}

\title{The Berwald-type linearisation of generalised connections}

\author{Tom Mestdag and Willy Sarlet\\
{\small Department of Mathematical Physics and Astronomy }\\
{\small Ghent University, Krijgslaan 281, B-9000 Ghent, Belgium}}
\date{}

\maketitle

\begin{quote}
{\bf Abstract.} {\small We study the existence of a natural
`linearisation' process for generalised connections on an affine
bundle. It is shown that this leads to an affine generalised
connection over a prolonged bundle, which is the analogue of what
is called a connection of Berwald type in the standard theory of
connections. Various new insights are being obtained in the fine
structure of affine bundles over an anchored vector bundle and
affineness of generalised connections on such bundles.}
\end{quote}


\section{Introduction}

The notion of Berwald connection seems to have its origin in Finsler
geometry. A Finsler spray generates a non-linear connection on
the tangent bundle $\tau_M:TM \rightarrow M$ and the Finslerian
Berwald connection represents a linearised version of this
non-linear connection. It is by now well known, however, that this
linearisation process can be applied to any non-linear connection
on the tangent bundle, resulting in a connection which is said to
be `of Berwald type'. There are different equivalent descriptions
in the literature concerning this process of linearisation, in
particular concerning the kind of bundle on which this is taking
place (see e.g.\ \cite{Szil,Bejancu,Crampin}). We adopt the line
of thinking of, for example, \cite{MarCar,Crampin}, where the linear
connection associated to a non-linear one on $\tau_M$, is regarded as
a connection on the pullback bundle $\tau_M^*TM\rightarrow TM$.

Much earlier already, Vilms showed \cite{Vilms} that
Berwald-type connections can be constructed starting from a
non-linear connection on an arbitrary vector bundle, not
necessarily a tangent bundle. In this paper, we wish to extend the
notion of Berwald connection even further by investigating the
following kind of generalisations. In the first place, motivated by
the recent interest in so-called `generalised connections', we shall
explore how the linearisation idea works in that context.
Generalised connections are connections on a bundle $\pi:E \rightarrow
M$ over a vector bundle ${\sf V}\rightarrow M$ say, which is anchored
in $TM$ via a bundle morphism $\rho:{\sf V}\rightarrow TM$ (think, for example,
of a Lie algebroid). Such connections show up in various fields of
application: see for example the work of Fernandes on Poisson geometry
\cite{Fernandes2}, and recent applications discussed by Langerock in the
fields of non-holonomic mechanics \cite{BavoI}, sub-Riemannian geometry
\cite{BavoII}, and control theory \cite{BavoIII}. For a general account
on generalised connections, see \cite{FransBavo} and
references therein. In many applications, the bundle $\pi$ will itself
be a vector bundle. The second kind of generalisation we wish to
investigate, however, is the situation where $\pi$ is an affine bundle.
In fact, this paper will almost entirely deal with the case of an affine
bundle, because there are lots of subtle points to be understood in such
a case, while it is easy to deduce the corresponding results for vector
bundles from the affine case.

An obvious motivation for paying attention to the case of affine bundles
comes from the geometry of time-dependent mechanical systems.
Time-dependent second-order ordinary differential equations, for example,
(\sode s for short) are modelled by a vector field on the first-jet bundle
of a manifold which is fibred over the real numbers: this is an
affine bundle over the base manifold. \sode s provide a canonically
defined non-linear connection on this affine bundle and the associated
linear connection has played a role in a variety of applications.
In \cite{MeSa}, we have made a quite exhaustive comparative study of
different versions of this linear connection (not always called a
Berwald-type connection in the literature), which were independently described
by Byrnes \cite{Byrnes}, Massa and Pagani \cite{MaPa}, and
Crampin {\em et al\/} \cite{CraMaSa}. As said before, our line of approach
was to view this linear connection as being defined on some pullback bundle.
We observed that the differences come from a kind of `gauge freedom' in fixing
the `time-component' of the connection. In particular, we found that there are
two rather natural ways of fixing this freedom. Our present
contribution is in some respect a continuation of this work. As a matter
of fact, within the much more general context of generalised connections
and arbitrary affine bundles, we will be led to a very clear understanding
of the origin of these two competing natural constructions. The present work
further links up with our recent studies of affine Lie algebroids and
`time-dependent Lagrangian systems' defined on such affine algebroids
\cite{SaMeMa,MaMeSa}. Last but not least, this paper complements (and
in fact was announced in) our recent analysis of affineness of generalised
connections \cite{MeSaMa}.

In Section 2, we recall the basic features of affine bundles and
generalised connections, needed for the rest of the paper.
Most of Section~3 is about the special case that the generalised connection
on the affine bundle is itself affine, in the sense of \cite{MeSaMa}.
We study the relationship between parallel transport along an admissible curve
and Lie transport of vertical vectors along the horizontal lift of such a
curve, and further arrive at explicit defining relations for the covariant
derivative operators associated to an affine generalised connection.
This paves the way to the analysis of the next section, where the idea
is to conceive a notion of Berwald-type connection associated to an
arbitrary (non-linear) generalised connection. This Berwald-type connection
indeed appears to be an affine generalised connection over the prolonged
bundle of the original affine bundle $\pi:E\rightarrow M$.
In the spirit of \cite{Crampin}, the way the idea of Berwald-type
connections is being developed, comes from looking at natural ways of
defining (on the pullback bundle $\pi^*E$) rules of parallel transport
along horizontal and vertical curves in $E$. There appear to be two
ways of defining a kind of notion of complete parallellism in the fibres
of $E$; they give rise to two different Berwald-type connections which
relate back to the results for time-dependent mechanics in \cite{MeSa}.
In the fifth section, we further specialise to the case
where the affine bundle is an affine Lie algebroid and the generalised
connection is the one canonically associated to a given pseudo-\sode,
as discussed already in \cite{Me}. The case of Lagrangian systems on
affine Lie algebroids then is a further particular situation.

\section{Basic setup}\label{basicsetupsection}

Let $\pi:E\rightarrow M$ be an affine bundle, modelled on the
vector bundle $\ov\pi:\ov E \rightarrow M$. The set of all affine
functions on $E_m$ ($m\in M$), $E^\dagger_m := \mbox{Aff}(E_m,\R)$
is the typical fibre of a vector bundle $E^\dagger= \bigcup_{m\in
M}E^\dagger_m$ over $M$, called the extended dual of $E$. In turn,
the dual of $\pi^\dagger: E^\dagger \rightarrow M$, denoted by
$\tpi: \tilde E:= (E^\dagger)^* \rightarrow M$, is a vector bundle
into which both $E$ and $\ov E$ can be mapped via canonical
injections. We will further need also the pullback bundles
$\pi^*\pi:E\rightarrow E$, $\pi^*\ov\pi:\ov E\rightarrow E$ and
$\pi^*\tpi:\tilde E \rightarrow E$. In the following, we will not
make a notational distinction between a point in $E$ and its
injection in $\tilde E$ (and likewise for a vector in $\ov E$).
Similarly, if for example $\sigma$ denotes a section of the affine
bundle, the same symbol will be used for its injection in
$Sec(\tpi)$ and even for the section $\sigma\circ\pi$ of
$\pi^*\tpi$ (if one looks at a section of $\pi^*\tpi$ as a map
$\tilde X:E\rightarrow \tilde E$ such that $\tpi\circ\tilde X =
\pi$). We trust that the meaning will be clear from the context.

The structure of the $\cinfty{E}$-module $Sec(\pi^*\tpi)$ deserves
some closer inspection. The injection of $E$ into $\tilde E$
provides a canonical section of $\pi^*\tpi$, which will be denoted
by $\cans$. Furthermore, there exists a canonical map $\vartheta:
Sec(\pi^*\tpi)\rightarrow Sec(\pi^*\tpi)$, which can be discovered
as follows (see \cite{MaMeSa}). First, within a fixed fibre
$\tilde{E}_m$, every $\tilde e$ defines a unique number,
$\lambda(\tilde e)$ say, determined by $\tilde e=\lambda(\tilde
e)e$ for some $e\in E_m$ ($\lambda(\tilde e)$ is zero when $\tilde
e$ belongs to $\ov{E}_m$ and one when $\tilde e\in E_m$). Thus,
choosing an arbitrary $a\in E_m$, we get a map
$\vartheta_a:\tilde{E}_m\rightarrow\tilde{E}_m, \tilde e\mapsto
\tilde e-\lambda(\tilde e)a$, which actually takes values in ${\ov
E}_m$, and therefore a map
\begin{equation}
\vartheta:\pi^*\tilde E\rightarrow \pi^*\ov E \subset \pi^*\tilde
E,\ (a,\tilde e)\mapsto (a,\tilde e-\lambda(\tilde e)a).
\label{canonicalmap}
\end{equation}
We will use the same notation for the extension of this map to
$Sec(\pi^*\tpi)$, i.e.\ for $\tilde X \in Sec(\pi^*\tpi)$,
$\vartheta(\tilde X)(e)=\vartheta(\tilde X(e))$. It follows that
every $\tilde X\in Sec(\pi^*\tpi)$ can be written in the form
\begin{equation}
\tilde X = f_{\tilde X}\,\cans + \vartheta(\tilde X), \qquad
\mbox{with} \quad f_{\tilde X}\in\cinfty{E}: \ f_{\tilde
X}(e)=\lambda(\tilde X(e)). \label{decomp}
\end{equation}
Clearly, if $\tilde X=\vartheta(\tilde X)$ in some open
neighbourhood in $E$, it means that $\lambda(\tilde X (e))=0$, so
that $\tilde X(e)\in\ov{E}$ in that neighbourhood, and such a
$\tilde X$ cannot at the same time be in the span of $\cans$. We
conclude that locally:
\begin{equation}
Sec(\pi^*\tpi) = \langle\cans\rangle \oplus Sec(\pi^*\ov\pi).
\label{decomp2}
\end{equation}
As a consequence, if $\{{\ov\sigma}_\alpha \}$ is a local basis
for $Sec(\ov\pi)$, then $\{\cans,{\ov\sigma}_\alpha \}$ is a local
basis for $Sec(\pi^*\tilde\pi)$.

For any affine bundle, there exists a well-defined notion of {\em
vertical lift} (cf.\ \cite{MaMeSa}). The vertical lift of a vector
${\ov e} \in \ov E_m$ to an element of $T_e E$ at a point $e \in
E_m$ is the point $v(e,\ov e)$ determined by the requirement that
for all functions $f\in \cinfty{E}$:
\[
v(e,\ov e)f = \frac{d}{dt}f(e+t{\ov e})|_{t=0}.
\]
Then, any $(e,\tilde e)\in\pi^*\tilde{E}$ can be vertically lifted
to the point $v(e,\tilde e)$ in the fibre $T_e E$ over $e$,
determined by
\[
v(e,\tilde e)=v(e,\vartheta_e(\tilde e)).
\]
The final step of course is to extend this construction in the
obvious way to an operation:
\[
v: Sec(\pi^*\tpi)\rightarrow \vectorfields E.
\]
It follows in particular that
\begin{equation}
v(\cans) = 0. \label{canvfv}
\end{equation}
Given a vertical $Q$ in $TE$, we will use the notation $Q_v$ for
the unique element in $\pi^*\ov E$ such that $v(Q_v)=Q$.

In what follows, there will be a role also for a second vector
bundle $\tau:{\sf V}\rightarrow M$ which is anchored in $TM$ by
means of a linear bundle map $\varrho: {\sf V} \rightarrow TM$.
Note that $\varrho$ can be regarded in an obvious way also as a
map (with the same symbol) from $\pi^*{\sf V}$ into $\pi^*TM$, by
means of $\varrho(e,{\sf v}) = (e,\varrho({\sf v}))$ for any
$(e,{\sf v})\in \pi^*{\sf V}$. The generalised connections we will
be concerned with in the rest of this paper are so-called {\em
$\varrho$-connections on $\pi$}. They are defined (see
\cite{FransBavo} and references therein) as follows.
\begin{dfn}
A $\varrho$-connection on $\pi$ is a linear bundle map
$h:\pi^*{\sf V}\rightarrow TE$ such that $T\pi\circ h= \varrho
\circ p_{\sf V}$.
\end{dfn}
Here, $p_{\sf V}$ is the projection $\pi^*{\sf V}\rightarrow {\sf
V}$. In \cite{MeSaMa} we have shown that the terminology
`connection' is justified here, since a $\varrho$-connection can
be seen, alternatively, as a splitting of some short exact
sequence. For that purpose, one has to invoke the {\em
$\varrho$-prolongation of $\pi$} (see e.g. \cite{Mart,Higgins}).
It is the bundle $\pi^1: \prolr{E} \rightarrow E$ whose total
space $\prolr{E}$ is the total space of the pullback bundle
$\varrho^*TE$,
\begin{equation}
\prolr{E} = \{({\sf v},Q_e) \in {\sf V} \times TE \mid \,
\varrho({\sf v}) = T\pi(Q_e) \} \label{LVP}
\end{equation}
whereby the projection $\pi^1$ is the composition of the
projection $\varrho^1$ of $\varrho^*TE$ onto $TE$ with the tangent
bundle projection $\tau_E$, $\pi^1=\tau_E\circ\varrho^1$. The
vector bundle $\pi^1$ has a well-defined subbundle ${\mathcal
V}^\varrho E \rightarrow E$, the vertical bundle, consisting of
those elements that lie in the kernel of the projection $\prolr
E\rightarrow {\sf V}$. These elements are of the form $(0,Q_e)$,
where $Q_e$ is also vertical in $T_eE$. We can now extend the
vertical lift $v$ to a map $\V{}:\pi^*\tilde E\rightarrow \prolr
E$, by means of
\begin{equation}
(e,\tilde e\V)=(0,v(e,\tilde e)). \label{vlift}
\end{equation} The point now is
that a $\varrho$-connection on $\pi$ can equivalently be seen as a
splitting of the short exact sequence
\begin{equation}
0 \rightarrow {\mathcal V}^\varrho E \stackrel{}{\rightarrow}
\prolr{E} \stackrel{j}{\rightarrow} \pi^* {\sf V} \rightarrow 0,
\label{shortexactgeneral}
\end{equation}
with $j: \prolr{E} \rightarrow \pi^* {\sf V}: ({\sf v},Q) \mapsto
(\tau_E(Q),{\sf v})$, i.e. as a map $\H{}:\pi^*{\sf V}\rightarrow
\prolr E$ such that $j\circ\H{}=id_{\pi^*{\sf V}}$. The relation
between the maps $h$ and $\H{}$ is: $\varrho^1\circ \H{}=h$. As
always, we will use the same symbol for the extension of the maps
$h$, $v$, $\H{}$ and $\V{}$ to sections of the corresponding
bundles. As a consequence of the existence of a splitting, for any
section ${\mathcal Z}\in Sec(\pi^1)$, there exist uniquely
determined sections ${\sf X}\in Sec(\pi^*\tau)$ and $\ov Y\in
Sec(\pi^*\ov\pi)$ such that
\begin{equation}
{\mathcal Z} = \H{\sf X} + \V{\ov Y}. \label{decomposition}
\end{equation}
In fact, if $\{{\sf s}_a\}$ is a local basis for $Sec(\tau)$ and
$\{{\ov \sigma}_\alpha\}$ a basis for $Sec(\ov \pi)$, and these
are interpreted as sections of $\pi^*{\sf V}\rightarrow E$ and
$\pi^*\ov E \rightarrow E$, respectively, then $\{\H{\sf
s}_a,\V{\ov \sigma}_\alpha\}$ provides a local basis for
$Sec(\pi^1)$.

A summary of most spaces and maps involved in the above
construction is presented in the following diagram.

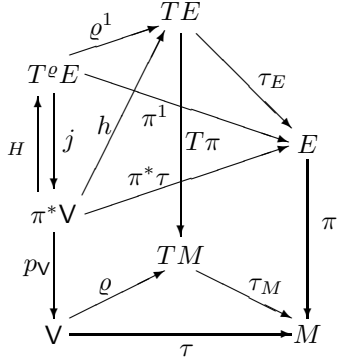
\setlength{\unitlength}{12pt}
\begin{picture}(40,13)(-2,2)
\put(7.5,7.3){\vector(3,1){6.4}}
\put(7.5,11.7){\vector(3,-1){6.4}} \put(11,13){\vector(1,-1){3}}
\put(7,12.3){\vector(3,1){2.8}} \put(7.4,7.9){\vector(1,2){2.6}}
\put(5.2,9.2){\mybox{ $^H$}}
 \put(6.0,8.0){\vector(0,1){3}}
 \put(6.5,7.3){\mybox{$\pi^*{\sf V}$}} \put(6.5,11.7){\mybox{$\prolr{E}$}}
\put(14.5,9.5){\mybox{$E$}} \put(10.5,13.7){\mybox{$TE$}}
\put(7,3.5){\vector(1,0){7}} \put(11,5.5){\vector(2,-1){3}}
\put(7,4){\vector(2,1){3}} \put(6.5,3.5){\mybox{${\sf V}$}}
\put(14.5,3.5){\mybox{$M$}} \put(10.5,6){\mybox{$TM$}}
\put(6.5,6.7){\vector(0,-1){2.7}}
\put(10.5,13){\vector(0,-1){6.4}} \put(14.5,9){\vector(0,-1){5}}
\put(6.5,11.1){\vector(0,-1){3.0}} \put(10.7,3){\mybox{$\tau$}}
\put(8.1,5.0){\mybox{$\varrho$}} \put(13.2,5.0){\mybox{$\tau_M$}}
\put(9.5,8.5){\mybox{$\pi^*\tau$}}
\put(8.0,13.2){\mybox{$\varrho^1$}}
\put(13.4,11.5){\mybox{$\tau_E$}} \put(9.7,10.4){\mybox{$\pi^1$}}
\put(6.0,5.6){\mybox{$p_{\sf V}$}} \put(15.2,7){\mybox{$\pi$}}
\put(7.0,9.6){\mybox{$j$}} \put(11.2,9.6){\mybox{$T\pi$}}
\put(8.1,10.1){\mybox{$h$}}
\end{picture}
\setlength{\unitlength}{\savelen}

For later use we list here some coordinate expressions. Let $x^i$
denote coordinates on $M$ and $y^\alpha$ fibre coordinates on $E$
with respect to some local frame $\cosyst$ for $Sec(\pi)$. For
each $\sigma\in Sec(\pi)$ with local representation
$\sigma(x)=e_0(x) + \sigma^\alpha(x){\ov e}_\alpha(x)$, the maps
\begin{equation}
e^0(\sigma)(x)=1,\ \forall x,\qquad
e^\alpha(\sigma)(x)=\sigma^\alpha(x).\label{basis}
\end{equation}
define an induced basis for $Sec(\pi^\dagger)$. Remark that $e^0$
coincides, in each fibre, with the constant function 1 and thus
has a global character. We will denote by $(e_0,e_\alpha)$ the
basis of $Sec(\tpi)$ dual to the basis (\ref{basis}). Since
sections of $\tpi$ can be regarded also as (basic) sections of
$\pi^*\tpi$, $(e_0,e_\alpha)$ can serve at the same time as local
basis for $Sec(\pi^*\tpi)$. Hence, every $\tilde X\in
Sec(\pi^*\tpi)$ can be represented in the form $\tilde X={\tilde
X}^0(x,y)e_0 + {\tilde X}^\alpha(x,y) e_\alpha$. But more
interestingly, with the use of the canonical section $\cans$, we
have
\begin{equation}
\cans = e_0 + y^\alpha e_\alpha, \qquad \tilde X={\tilde X}^0\cans
+ ({\tilde X}^\alpha-{\tilde X}^0y^\alpha)e_\alpha.
\label{candecomp}
\end{equation}

If we use ${\sf v}^a$ for the fibre coordinates of $\tau$ and $\{
{\sf e}_a \}$ for the corresponding local basis, then a given
anchor map $\varrho$ has coordinate representation $\varrho:
(x^i,{\sf v}^a) \mapsto \varrho^i_a(x) {\sf v}^a\fpd{}{x^i}$. In
\cite{MeSaMa} we have shown that a natural choice for a local
basis of sections of the prolonged bundle
$\pi^1:\prolr{E}\rightarrow E$ is given by the following: for each
$e\in E$, if $x$ are the coordinates of $\pi(e)\in M$,
\begin{equation}
\baseX a (e)= \left({\sf e}_a(x),\left.\varrho^i_a(x)\fpd{}{x^i}
\right|_e\right), \quad \baseV \alpha(e) =
\left(0,\left.\fpd{}{y^{\alpha}}\right|_e\right).
\label{prolongedbasis}
\end{equation}
The dual basis for $Sec({\pi^1}{}^*)$ is denoted by $\{ {\mathcal
X}^a, {\mathcal V}^\alpha \}$. A general section of the prolonged
bundle can be represented locally in the form:
\begin{equation}
{\mathcal Z}= {\sf z}^a(x,y)\baseX{a} +
Z^\alpha(x,y)\baseV{\alpha}. \label{prolZ}
\end{equation}

Suppose that, in addition, we have a $\varrho$-connection on $\pi$
at our disposal. As in \cite{MeSaMa}, the local expressions of $h$
and $\H{}$ then are:
\[
h(x^i,y^\alpha,{\sf v}^a)=(x^i,y^\alpha,\rho^i_a(x){\sf
v}^a,-\Gamma^\alpha_a(x,y){\sf v}^a),
\]
and
\[
\H{(x^i,y^\alpha,{\sf v}^a)}= \left((x^i,{\sf v}^a),{\sf
v}^a\left(\rho^i_a\fpd{}{x^i}
-\Gamma^\alpha_a\fpd{}{y^\alpha}\right)\right).
\]
We can now easily give a local basis for the horizontal sections
of $\pi^1$, which is given by
\begin{equation}
\baseH{a}=\H{{\sf e}_a}=\baseX{a} - \Gamma_a^\alpha(x,y) \baseV
\alpha. \label{hor}
\end{equation}
An adapted representation of the section (\ref{prolZ}) then
becomes:
\begin{equation}
{\mathcal Z}=  {\sf z}^a\baseH{a} + (Z^\alpha+{\sf
z}^a\Gamma^\alpha_a) \baseV{\alpha}. \label{prolZ2}
\end{equation}

We end this section with the following remark about a generalised
notion of tangent map between $\varrho$-prolongations. Suppose
that two vector bundles ${\sf V}_1 \rightarrow M_1$ and ${\sf V}_2
\rightarrow M_2$ with anchors $\varrho_1$ and $\varrho_2$
(respectively) are given, together with two arbitrary fibre
bundles $P_1\rightarrow M_1$ and $P_2 \rightarrow M_2$. Suppose
further that $F: P_1 \rightarrow P_2$ is a bundle map over some
$f:M_1 \rightarrow M_2$ and that ${\sf f}: {\sf V}_1 \rightarrow
{\sf V}_2$ is a vector bundle morphism over the same $f$,
satisfying $Tf\circ\varrho_1 = \varrho_2\circ {\sf f}$. Then, we
can define a map
\begin{equation}
{T}^{\varrho_1,\varrho_2}F : T^{\varrho_1} P_1 \rightarrow
{T}^{\varrho_2} P_2, ({\sf v}_1,X_1) \mapsto ({\sf f}({\sf
v}_1),TF(X_1)). \label{LV1V2F}
\end{equation}

\section{Parallel transport and Lie transport for affine generalised connections}

Let $h: \pi^*{\sf V}\rightarrow TE$ be a $\varrho$-connection on
an affine bundle $\pi: E \rightarrow M$. The affine structure of
$E$ can be represented by the map $\Sigma: E \times_M E
\rightarrow E$, $\Sigma(e,\ov e)=e+\ov e$ ($e\in E_m$, $\ov e\in
{\ov E}_m$). The $\varrho$-connection $h$ on $\pi$ is said to be
{\em affine}, if there exists a linear $\varrho$-connection $\ov
h$ on $\ov\pi$ such that $\forall e\in E_m$, $\ov e \in {\ov E}_m$
and ${\sf v} \in {\sf V}_m$
\begin{equation}
h(e+\ov e,{\sf v}) = T\Sigma_{(e,\ov e)}\Big(h(e,{\sf v}),\ov h
(\ov e, {\sf v})\Big). \label{affineconnection}
\end{equation}
The connection coefficients of an affine connection are of the
form $\Gamma_a^\alpha(x,y) = \Gamma_{a0}^{\alpha}(x)+
\Gamma_{a\beta}^{\alpha}(x)y^\beta$. We have shown in
\cite{MeSaMa} that an affine $\varrho$-connection on $\pi$ can
equivalently be seen as a pair $(\nabla,\ov\nabla)$ where
$\ov\nabla$ is the covariant derivative operator corresponding to
the linear $\ov h$ and $\nabla: Sec(\tau)\times
Sec(\pi)\rightarrow Sec(\ov\pi)$ is an operator which is
$\R$-linear in its first argument and has the properties
\begin{equation}
\nabla_{f{\sf s}}\sigma = f  \nabla_{{\sf s}}\sigma, \qquad
\nabla_{\sf s}(\sigma + f\ov \sigma)=\nabla_{{\sf s}}\sigma +
f\ov\nabla_{{\sf s}}\ov\sigma + \varrho({\sf s})(f)\ov\sigma
\end{equation}
for all ${\sf s}\in Sec(\tau)$, $\sigma\in Sec(\pi)$, $\ov\sigma
\in Sec(\ov\pi)$ and $f\in \cinfty M$. $\nabla$ is defined via the
so-called connection map $K: \prolr E \rightarrow \ov E$, which
maps $({\sf v}, Q_e)$ to the vertical tangent vector
$(\varrho^1-h\circ j)({\sf v, Q_e})$ thought of as an element of
$\ov E$,
\[
\nabla_{\sf s}\sigma(m)=K\big({\sf s}(m), T\sigma(\varrho({\sf
s}(m)))\big).
\]
Technically, $K=p_{\ov E}\circ {}_v\circ(\varrho^1-h\circ j)$.

Let us briefly summarise the discussion of parallel transport for
affine $\varrho$-connections, as developed in \cite{MeSaMa}. A
curve in $\pi^*{\sf V}$ is a couple $(\psi,{\sf c})$, where $\psi$
is a curve in $E$ and ${\sf c}$ a curve in ${\sf V}$ with the
properties that the projected curves on $M$ coincide: $\psiM=\cM$
(in taking a curve in $\pi^*{\sf V}$ we will suppose that
$I=[a,b]\subset Dom((\psi,{\sf c}))$ is an interval in
$Dom(\psi)\cap Dom({\sf c})$). The {\em horizontal lift of the
curve $(\psi,{\sf c})$} is a curve $(\psi,{\sf c}\H)$ in $\prolr
E$, determined by
\begin{equation}
(\psi,{\sf c}\H): u \mapsto ({\sf c}(u),h(\psi(u),{\sf c}(u)))
\qquad \mbox{for all $u\in I$.} \label{star}
\end{equation}
We also say that $\psi$ in $E$ is {\em a horizontal lift of ${\sf
c}$}, and write $\psi=c^h$, if $(\psi,{\sf c}\H)$ is a
$\varrho^1$-admissible curve. Since by construction
$\pi^1\circ(\psi,{\sf c}\H)=\psi$, this means that
\begin{equation} \label{dotpsi}
\dot\psi(u) = \varrho^1\circ(\psi,{\sf c}\H)(u)=h(\psi(u),{\sf
c}(u)),
\end{equation}
and automatically implies that ${\sf c}$ must be
$\varrho$-admissible, since $\dot{c}_{\scriptscriptstyle
M}=\dot{\psi}_{\scriptscriptstyle M}=T\pi\circ \dot\psi = T\pi
\circ h(\psi,{\sf c})= \varrho \circ {\sf c}$. Given ${\sf c}$,
with $\cM(a)=m$ say, we will write $c^h_e$ for the unique solution
of (\ref{dotpsi}) passing through the point $e\in E_m$ at $u=a$
(i.e. $c^h_e(a)=e$) and denote the lift $(c^h_e,{\sf c}\H)$ by
$\H{{\dot c}_e}$, where, of course, the `dot' merely refers to the
fact that $\varrho^1(\H{{\dot c}_e})={\dot c}^h_e$. Curves of the
form $\H{{\dot c}_e}$ are $\varrho^1$-admissible by construction.

Let ${\sf s}$ be a section of ${\sf V}\rightarrow M$, which we
regard as section of $\pi^*{\sf V}\rightarrow E$ via the
composition with $\pi$. In that sense, we can talk about the
vector field $h({\sf s})\in \vectorfields E$, which has the
following interesting characteristics.
\begin{lemma} \label{Lemma1}
For any ${\sf s} \in Sec(\tau)$, the vector field $h({\sf s})$ on
$E$ has the property that all its integral curves are horizontal
lifts of $\varrho$-admissible curves in ${\sf V}$.
\end{lemma}
\proof\, Let $\gamma_e$ denote an integral curve of $h(\sf s)$
through the point $e$. Since $h({\sf s})$ is $\pi$-related to
$\varrho({\sf s})\in \vectorfields M$, $\pi\circ\gamma_e$ then is
an integral curve of $\varrho({\sf s})$ through $m=\pi(e)$, which
we shall call $c_m$. Obviously, ${\sf c}={\sf s}(c_m)$ now is a
$\varrho$-admissible curve in ${\sf V}$ and we have for all $u$ in
its domain,
\begin{equation}
h({\sf s})(c^h_e(u))=h\Big(c^h_e(u),{\sf
s}\big(\pi(c^h_e(u))\big)\Big) = h\big(c^h_e(u), {\sf c}(u)\big) =
{\dot c}^h_e(u), \label{fourtytwo}
\end{equation}
which shows that $\gamma_e=c^h_e$. \qed

So far, the above characterisation of a horizontal lift applies to
any $\varrho$-connection on $\pi$. If, in particular, the
connection is {\em affine}, then we can express the definition of
a horizontal lift also in terms of the operator $\nabla$. Indeed,
for any $\varrho$-admissible curve ${\sf c}$ in ${\sf V}$ and any
curve $\psi$ with $\cM=\psiM$, we can define a new curve
$\nabla_{\sf c}\psi$ by means of
\[
\nabla_{\sf c}\psi(u):=K(({\sf c}(u),\dot\psi(u))) = \Big(\dot\psi
(u) - h\big((\psi(u),{\sf c}(u))\big)\Big)_v.
\]
Obviously, $\psi=c^h_e$ iff $\psi(a)=e$ and $\nabla_{\sf c}\psi=0$
for all $u\in I$.

Putting $\cM(b)=m_b$, the point $c^h_e(b)\in E_{m_b}$ is called
the {\em parallel translate of e along ${\sf c}$.} It is
instructive to see in detail how one can get an affine action
between the affine fibres of $E$. Take $e_1,e_2\in E_{m}$ and
consider the horizontal lifts $c^h_{e_1}$ and $c^h_{e_2}$. Denote
the difference $e_1-e_2$ by $\ov e \in \ov E_{m}$ and put
$\ov\eta_{\ov e}:= c^h_{e_1}-c^h_{e_2}$. As the subscript
indicates, $\ov\eta_{\ov e}$ is a curve in $\ov E$ starting at
$\ov e$. From the action of $\nabla$ on curves (see e.g.
\cite{FransBavo,MeSaMa}), it easily follows that
\[
\nabla_{\sf c} c^h_{e_1}(u)- \nabla_{\sf c} c^h_{e_2}(u) =
{\ov\nabla}_{\sf c} {\ov\eta}_{\ov e}(u), \qquad \mbox{for all
$u\in I$.}
\]
Since both $c^h_{e_1}$ and $c^h_{e_2}$ are solutions of the
equation $\nabla_{\sf c} \psi=0$, ${\ov\eta}_{\ov e}$ must be the
unique solution of the initial value problem $\ov\nabla_{\sf
c}\ov\eta =0$, $\ov\eta(a)=\ov e$, i.e. the unique $\ov
h$-horizontal lift through $\ov e$. Therefore, {\em the difference
between the $\nabla$-parallel translates of $e_1$ and $e_2$ along
${\sf c}$ is the $\ov\nabla$-parallel translated of $e_1-e_2$
along ${\sf c}$}. In fact, this property is necessary and
sufficient for the connection to be {\em affine}. From now on we
will use the notation $c^{\ov h}_{\ov e}$ for
$c^h_{e_1}-c^h_{e_2}$.
\begin{prop}
A $\varrho$-connection $h$ on $\pi$ is affine if and only if there
exists a linear $\varrho$-connection $\ov h$ on $\ov\pi$, such
that for all admissible curves ${\sf c}$ and for any two points
$e_1,e_2\in E_m$ $(m=\cM(a))$ the difference
$c^h_{e_1}(u)-c^h_{e_2}(u)$ is the $\ov\nabla$-parallel translate
of $e_1-e_2$ along ${\sf c}$. \label{propone}
\end{prop}
\proof\, The proof in one direction has already been given. For
the converse, suppose that a {\em linear} connection $\ov h$
exists, having the above properties. It suffices to show that $\ov
h$ is related to $h$ by means of (\ref{affineconnection}).
Choosing a ${\sf v}\in {\sf V}_m$ arbitrarily, we take a
$\varrho$-admissible curve ${\sf c}$, that passes through it (for
$u=a$) and consider its $h$-horizontal lift $c^h_e$ through $e$
and its $\ov h$-horizontal lift through $\ov e$. By assumption we
know that $c^h_{e+\ov e}- c^h_{e}$ is the $\ov h$-horizontal lift
of ${\sf c}$ through $\ov e$, i.e.
\[
\Sigma(c^h_e(u),{c}^{\ov h}_{\ov e}(u))=c^h_{e+\ov e}(u) \qquad
\mbox{for all $u\in I$.}
\]
Taking the derivative of this expression at $u=a$, we get
\begin{equation}
T\Sigma_{(e,\ov e)}\big({\dot c}^h_e(a), {\dot {c}}^{\ov h}_{\ov
e}(a)\big) = {\dot c}^h_{e+\ov e}(a). \label{Tsigmachovch}
\end{equation}
In view of (\ref{dotpsi}) and its analogue for $\ov h$, this is
indeed what we wanted to show. \qed

Let us introduce corresponding `flow-type' maps. For that purpose,
it is convenient to use (temporarily) the more accurate notation
$c^h_{a,e}$ for the horizontal lift which passes through $e$ at
$u=a$. Putting
\[
\phi^h_{u,a}(e)=c^h_{a,e}(u), \qquad u\in [a,b],
\]
the result of the preceding proposition can equivalently be
expressed as,
\[
\phi^h_{u_2,u_1} (e+\ov e) = \phi^h_{u_2,u_1}(e)+  \phi^{\ov
h}_{u_2,u_1}(\ov e),
\]
i.e. $\phi^h_{u_2,u_1}: E_{\cM(u_1)}\rightarrow E_{\cM(u_2)}$ is
an affine map with linear part $\phi^{\ov h}_{u_2,u_1}$. Its
tangent map, therefore, can be identified with its linear part. As
a result, when we consider Lie transport of vertical vectors in
$TE$ along the curve $c^h_{a,e}$, in the case of an affine
connection, the image vectors will come from the parallel
translate associated to the linear connection $\ov h$. Indeed,
putting $Y_a=v(e,\ov e)$ and defining its Lie translate from $a$
to $b$ as $Y_b=T\phi^h_{b,a}(Y_a)$, we have for each $g\in \cinfty
E$,
\begin{eqnarray*}
Y_b(g)= Y_a(g \circ \phi^h_{b,a}) &=& \frac{d}{dt}\big(g\circ
\phi^h_{b,a}(e+t\ov e)\big)_{t=0} \\ &=&
\frac{d}{dt}\big(g(\phi^h_{b,a}(e)+t\phi^{\ov h}_{b,a}(\ov
e))\big)_{t=0} = v\big(c^h_{a,e}(b),c^{\ov h}_{a,\ov e}(b)\big)g,
\end{eqnarray*}
where the transition to the last line requires affineness of the
connection. It follows that in the affine case, the Lie translate
of $Y_a=v(e,\ov e)$ to $b$ is given by
\begin{equation}
Y_b=v\big(c^h_{a,e}(b),c^{\ov h}_{a,\ov e}(b)\big).
\label{starstar}
\end{equation}

At this point, it is appropriate to make a few more comments about
the general idea of Lie transport. If (on an arbitrary manifold)
$Y$ is a vector field along an integral curve of some other vector
field $X$, and we therefore have a genuine (local) flow $\phi_s$
at our disposal, then the Lie derivative of $Y$ with respect to
$X$ is defined to be
\[
{\mathcal L}_X Y (u) =
\frac{d}{ds}\big(T\phi_{-s}(Y(s+u))\big)_{s=0}.
\]
As shown for example in \cite{CrampinPirani} (p.\ 68), the
requirement ${\mathcal L}_X Y=0$, subject to some initial
condition, $Y(0)=Y_0$ say, then uniquely determines a vector field
$Y$ along an integral curve of $X$ in such a way that $Y(u)$ is
obtained by Lie transport of $Y_0$. The description of Lie
transport, therefore, is more direct when we are in the situation
of an integral curve of a vector field.

Lemma~\ref{Lemma1}, unfortunately, does not create such a
situation for us because, when an admissible curve ${\sf c}$ is
given, together with one of its horizontal lifts $c^h_e$, it does
not provide us with a way of constructing a vector field which has
$c^h_e$ as one of its integral curves. The complication for
constructing such a vector field primarily comes from the fact
that the differential equations (\ref{dotpsi}) which define
$c^h_e$ are non-autonomous. The usual way to get around this
problem is to make the system autonomous by adding on extra
dimension. A geometrical way of achieving this here, which takes
into account that $c^h_e$ in the first place has to be a curve
projecting onto $\cM$, is obtained by passing to the pullback
bundle $\cM^*E \rightarrow I \subset \R$. We introduce the
notation $\cM^1:\cM^* E \rightarrow E, (u,e)\mapsto (\cM(u),e)$,
and likewise ${\ov c}_{\scriptscriptstyle M}^1:\cM^* \ov E
\rightarrow \ov E, (u,\ov e)\mapsto (\cM(u),\ov e)$. With the help
of these maps, we can single out vector fields $\Lambda_c \in
\vectorfields{\cM^*E}$ and ${\ov\Lambda}_c\in
\vectorfields{\cM^*\ov E}$ as follows.
\begin{prop}
For any $\varrho$-connection $h$ on $\pi$ and given
$\varrho$-admissible curve ${\sf c}$ in ${\sf V}$, there exists a
unique vector field $\Lambda_c$ on $\cM^*E$ that projects on the
coordinate vector field on $\R$ and is such that
\begin{equation}
T\cM^1(\Lambda_c(u,e))= h\big(\cM^1(u,e),{\sf
c}(u)\big),\label{Tc1MSigmac}
\end{equation}
for all $(u,e) \in \cM^*E$. Likewise, if the connection is affine
with linear part $\ov h$, there exists a unique vector field
${\ov\Lambda}_c$ on $\cM^*{\ov E}$ that projects on the coordinate
vector field on $\R$ and is such that
\begin{equation}
T{\ov c}_{\scriptscriptstyle M}^1({\ov \Lambda}_c(u,{\ov e}))=
{\ov h}\big({\ov c}_{\scriptscriptstyle M}^1(u,\ov e),{\sf
c}(u)\big),\label{Tc1MSigmac2}
\end{equation}
for all $(u,\ov e) \in \cM^*\ov E$. \label{Sigmaprop}
\end{prop}
\proof \, The proof is analogous for both cases; we prove only the
first. Let $u$ denote the coordinate on $\R$ and $y^\alpha$ the
fibre coordinates of some $e\in (\cM^*E)_u$. Representing the
given curve as ${\sf c}: u \mapsto (x^i(u),{\sf c}^a(u))$ and
putting
$\Lambda_c(u,e)=U(u,e)\fpd{}{u}|_{(u,e)}+Y^\alpha(u,e)\fpd{}{y^\alpha}|_{(u,e)}$,
one finds that
\[
T\cM^1(\Lambda_c(u,e)) = (x^i(u),y^\alpha; U(u,e) {\dot x}^{i}(u),
Y^\alpha(u,e)).
\]
On the other hand $h\big(\cM^1(u,e),{\sf c}(u)\big) =
\big(x^i(u),y^\alpha; {\sf c}^a(u) \varrho^i_a(\cM(u)), -{\sf
c}^a(u) \Gamma_a^\alpha(\cM(u),e)\big)$. Identification of the two
expressions gives that $Y^\alpha(u,e)=-{\sf c}^a(u)
\Gamma_a^\alpha(\cM(u),e)$, and that $U{\dot x}^i=\varrho^i_a {\sf
c}^a$. At points where ${\sf c}$ does not lie in the kernel of
$\varrho$, the latter equality would by itself determine $U$ to be
1. The extra projectability requirement ensures that this will
hold also when ${\dot x}^i=\varrho^i_a{\sf c}^a=0$. \qed

An interesting point here is that the complication about ensuring
separately that $U=1$ in some sense disappears when we look at
horizontality on $\prolr E$ rather than on $TE$, i.e.
horizontality in the sense of (\ref{star}). To see this, observe
first that we can use the bundle map  ${\sf f}:T\R\rightarrow {\sf
V}, U\frac{d}{du}|_u \mapsto U {\sf c}(u)$ over $\cM: I
\rightarrow M$, to obtain, in accordance with (\ref{LV1V2F}), the
following extended notion of tangent map:
\[
T^{\varrho}\cM^1: T(\cM^* E) \rightarrow \prolr E, \lambda \mapsto
({\sf f}(T\tau_\R \lambda),T\cM^1(\lambda)).
\]
Here $\tau_\R$ is the bundle projection of $\cM^*E\rightarrow\R$.
$T^{\varrho}\cM^1$ is well defined, since $T\tau_\R\lambda$ is of
the form $U\frac{d}{du}|_u$ and $T\pi\circ T\cM^1(\lambda) =
T(\pi\circ\cM^1)(\lambda)= T(\cM\circ\tau_\R)(\lambda)=T\cM\circ
T\tau_\R(\lambda) = T\cM (U\frac{d}{du}|_u) = U {\dot
c}_{\scriptscriptstyle M}(u) = U \varrho({\sf c}(u)) =
\varrho({\sf f}(U\frac{d}{du}|_u))=\varrho({\sf
f}(T\tau_\R\lambda))$. Notice that this remains true also at
points where ${\sf c}$ lies in the kernel of $\varrho$. Now,
$\Lambda_c$ can be defined as the unique vector field on $\cM^*E$
for which
\[
T^\varrho \cM^1(\Lambda_c(u,e))= \big(\cM^1(u,e),{\sf
c}(u)\H{\big)} \qquad \mbox{for all $(u,e)\in \cM^*E$}.
\]
Indeed, the second component of this equality is just the
condition (\ref{Tc1MSigmac}) again, whereas the first component
says that $U{\sf c}(u)={\sf c}(u)$ and thus implies $U=1$.

We will look now at the integral curves of $\Lambda_c$ and
${\ov\Lambda}_c$. Since $\Lambda_c$ projects on the coordinate
field on $\R$, the integral curves are essentially sections of
$\cM^*E \rightarrow \R$. Let $\gamma_e$ denote the integral curve
going through $(a,e)$ at time $u=a$, so that
\begin{equation}
{\dot\gamma}_e(u)=\Lambda_c(\gamma_e(u)) \qquad \mbox{$\forall
u\in I'$}, \label{intcurvcond}
\end{equation}
where $I'$ is some interval, possibly smaller than the domain $I$
of ${\sf c}$. In a similar way, we will write ${\ov\gamma}_{\ov
e}$ for the integral curve of ${\ov\Lambda}_c$ through $(a,\ov
e)$.
\begin{prop}
For any $\varrho$-connection on $\pi$, the curve
$\cM^1\circ\gamma_{e}$ is the $h$-horizontal lift of ${\sf c}$
through $e$. Likewise, if the connection is affine with linear
part $\ov h$, the curve ${\ov c}_{\scriptscriptstyle
M}^1\circ{\ov\gamma}_{\ov e}$ is the $\ov h$-horizontal lift of
${\sf c}$ through $\ov e$. \label{cM1circgamma}
\end{prop}
\proof\, Again, we will prove only the first statement.
$\cM^1\circ\gamma_{e}$ is a curve in $E$ projecting on the curve
$\cM$ in $M$. Using (\ref{intcurvcond}) and (\ref{Tc1MSigmac}) we
find:
\[
\frac{d}{du}(\cM^1\circ\gamma_{e})(u) = T \cM^1({\dot\gamma}_e(u))
= T\cM^1(\Lambda_c(\gamma_e(u)))=h(\cM^1\circ\gamma_{e}(u),c(u)),
\]
which shows that $\cM^1\circ\gamma_e=c^h_e$. \qed

We now proceed to look at Lie transport along integral curves of
the vector field $\Lambda_c$ on $\cM^*E$. It is clear that
$\cM^*E\rightarrow\R$ is an affine bundle modelled on the vector
bundle $\cM^*\ov E \rightarrow \R$. As we know from Section~2
then, there is a vertical lift map, which we will denote by $\vl$
which maps elements of $(\cM^*E \times_\R \cM^*\ov E)$ to vertical
vectors of $T(\cM^*E)$. We will consider Lie transport of such
vertical vectors along integral curves of $\Lambda_c$. Starting
from $e\in E_{\cM(a)}$, $\ov E \in {\ov E}_{\cM(a)}$ and putting
$\Upsilon_{e,\ov e}(a)=\vl \big((a,e),(a,\ov e)\big)$, we know
that the condition ${\mathcal L}_{\Lambda_c}\Upsilon_{e,\ov e}=0$
uniquely defines a vector field along the integral curve
$\gamma_{(a,e)}$ of $\Lambda_c$ ($\gamma_{(a,e)}(a)=(a,e)$) which
takes the (vertical) value $\Upsilon_{e,\ov e}(a)$ at the point
$(a,e)$. As said before, the value of $\Upsilon_{e,\ov e}$ at any
later $u$ is the Lie translate of $\Upsilon_{e,\ov e}(a)$, that is
to say, we have $\Upsilon_{e,\ov
e}(u)=T\phi_{u,a}\big(\Upsilon_{e,\ov e}(a)\big)$, where
$\phi_{u,a}$ refers to the flow of $\Lambda_c$ ($\phi_{a,a}$ is
the identity). It is interesting to observe here that
$\Upsilon_{e,\ov e}$ is directly related to the Lie translate we
discussed before, of vertical vectors on $E$ along the horizontal
lift $c^h_{a,e}$. To be precise, with $Y_{e,\ov
e}(a)=T\cM^1\big(\Upsilon_{e,\ov e}(a)\big)=v(e,\ov e)$ and
defining $Y_{e,\ov e}(u)$ to be $T\phi^h_{u,a}\big(Y_{e,\ov
e}(a)\big)$ as before, we have (at any later time $u$ in the
domain of $\gamma_{(a,e)}$)
\begin{equation}
T\cM^1\big(\Upsilon_{e,\ov e}(u)\big)=Y_{e,\ov e}(u).
\label{fifteen}
\end{equation}
Indeed, it follows from Proposition \ref{cM1circgamma} (using here
again the somewhat more accurate notations which take the
``initial time'' $a$ into account), that
$\cM^1\circ\phi_{u,a}=\phi^h_{u,a}\circ\cM^1$. Therefore, we have
\begin{eqnarray*}
T\cM^1\big(\Upsilon_{e,\ov e}(u)\big) &=&
T\cM^1\big(T\phi_{u,a}(\Upsilon_{e,\ov e}(a))\big) =
T\phi^h_{u,a}\big(T\cM^1(\Upsilon_{e,\ov e}(a))\big) \\ &=&
T\phi^h_{u,a}\big(Y_{e,\ov e}(a)\big)=Y_{e,\ov e}(u).
\end{eqnarray*}
The case of an {\em affine} $\varrho$-connection is of special
interest. The affine nature of the maps $c^h_e(u)$, for fixed $u$,
implies via Proposition \ref{cM1circgamma} that the flow maps of
$\Lambda_c$ are also affine, or expressed differently that:
\begin{equation}
\gamma_{e+\ov e}(u)=\gamma_e(u)+{\ov\gamma}_{\ov e}(u).
\label{sixteen}
\end{equation}
We have shown already that in the affine case: $Y_{e,\ov
e}(u)=v\big(c^h_e(u),c^{\ov h}_{\ov e}(u)\big)$. The translation
of this result via the relation (\ref{fifteen}) means that we have
\begin{equation}
\Upsilon_{e,\ov e}(u)=\vl\big(\gamma_e(u),{\ov\gamma}_{\ov
e}(u)\big).
\end{equation}
Summarising the more interesting aspects of what we have observed
above, we can make the following statement
\begin{prop} \label{Prop4}
For an arbitrary $\varrho$-connection on $\pi$, Lie transport of
vertical vectors on $E$ along the horizontal lift $c^h_e$ is
equivalent to Lie transport of vertical vectors on $\cM^*E$ along
integral curves of the vector field $\Lambda_c$. In the particular
case that the connection is affine, the translates in both cases
are compatible with the affine nature of the flow maps between
fixed fibres.
\end{prop}

The next result concerns an important property of the covariant
derivative operators $\nabla$ and $\ov\nabla$ which become
available when the $\varrho$-connection is affine. The preceding
considerations about vector fields $\Lambda_c$ on $\cM^*E$ will
help to prove it in a purely geometrical way. The reader may wish
to skip this rather technical proof and pass to the remark
immediately following it.

\begin{prop} \label{prop5}
Let $h$ be an affine $\varrho$-connection. For all ${\sf s}\in
Sec(\tau)$, $\ov\sigma\in Sec(\ov\pi)$ and $\sigma\in Sec(\pi)$,
the brackets $[h{\sf s},v\ov\sigma]$ and $[h{\sf s},v\sigma]$ of
vector fields on $E$ are vertical and we have
\begin{equation}
\ov\nabla_{\sf s}\ov\sigma = [h{\sf s},v\ov\sigma]_v , \qquad
\nabla_{\sf s}\sigma=[h{\sf s},v\sigma]_v. \label{mm}
\end{equation}
\end{prop}
\proof\, We start with the bracket $[h{\sf s},v\ov\sigma]_v$.
Since the vector fields under consideration are $\pi$-related to
$\varrho({\sf s})$ and the zero vector field on $M$, respectively,
their Lie bracket is $\pi$-related to $[\varrho({\sf s}),0]=0$ and
is therefore vertical. If we project it down to $\ov E$, strictly
speaking by taking $([h{\sf s},v\ov\sigma](e))_v \in E\times_M\ov
E$ with $e\in E_m$ and looking at the second component, we obtain
an element of ${\ov E}_m$ which does not depend on the fibre
coordinates of $e$ (as we can see instantaneously by thinking of
coordinate expressions). In other words, $[h{\sf s},v\ov\sigma]_v$
gives rise to a section of $\ov\pi$ which has the following
properties: for all $f\in \cinfty M$
\[
[h(f{\sf s}),v\ov\sigma]_v = f [h{\sf s},v\ov\sigma]_v, \qquad
[h{\sf s},v(f\ov\sigma)]_v = f[h{\sf s},v\ov\sigma]_v +
\varrho({\sf s})(f)\ov\sigma.
\]
These are precisely the properties of the covariant derivative
operator $\ov\nabla_{\sf s}\ov\sigma$, from which it follows that
$L({\sf s},\ov\sigma)=\ov\nabla_{\sf s}\ov\sigma - [h{\sf
s},v\ov\sigma]_v$ is tensorial in ${\sf s}$ and $\ov\sigma$.

To prove that $L$ is actually zero, we will use a rather subtle
argument, which is based on the following considerations of a
quite general nature. If $L({\sf s},\ov\sigma)$ is tensorial, so
that for all $m$,  $L({\sf s},\ov\sigma)(m)$ depends on ${\sf
s}(m)$ and $\ov\sigma(m)$ only, then for each curve $\cM$ in $M$,
there exists a corresponding operator $l({\sf r},\ov\eta)$, acting
on arbitrary sections along the curve $\cM$, which is completely
determined by the property: if ${\sf s}\in Sec(\tau)$,
$\ov\sigma\in Sec(\ov\pi)$ and we put ${\sf r}={\sf
s}\vert_{\cM}$, $\ov\eta=\ov\sigma\vert_{\cM}$, then $l({\sf
r},\ov\eta)(u) = L({\sf s},\ov\sigma)(\cM(u))$. In turn, the value
of $L({\sf s},\ov\sigma)(m)$ at an arbitrary point $m$ can be
computed by choosing an arbitrary curve $\cM$ through $m$
($\cM(0)=0$ say), selecting sections ${\sf r}$ and $\ov\eta$ along
$\cM$ for which ${\sf r}(0)={\sf s}(m)$ and $\ov\eta(0)=\ov e =
\ov\sigma(m)$, and then computing $l({\sf r},\ov\eta)(0)$.

We apply this general idea in the following way. Starting from an
arbitrary $e\in E_m$, we know from Lemma~\ref{Lemma1} that the
integral curve of $h({\sf s})$ through $e$ is the horizontal lift
$c^h_e$ of some admissible curve ${\sf c}$ through ${\sf s}(m)$
(here the ``initial time'' $a$ is taken to be zero). Take ${\sf
r}$ to be this curve ${\sf c}$ (with projection $\cM$) and choose
$\ov\eta$ to be the curve $c^{\ov h}_{\ov e}$. Then,
\[
l({\sf r},\ov\eta)(0)= \ov\nabla_{\sf c}c^{\ov h}_{\ov e}(0)-
\big( {\mathcal L}_{h{\sf s}} v(c^h_e,c^{\ov h}_{\ov e}) \big)_v
(0).
\]
But $\ov\nabla_{\sf c}c^{\ov h}_{\ov e}$ is zero by construction.
Concerning the second term, we observe that $h{\sf s}$ is
$\cM^1$-related to $\Lambda_c$, by definition of $\Lambda_c$. Also
$v(c^h_e,c^{\ov h}_{\ov e})$ is $\cM^1$-related to the vector
field $\Upsilon_{e,\ov e}(u)$ along the integral curve $\gamma_e$
of $\Lambda_c$ through the point $(0,e)\in \cM^*E$ (see
(\ref{fifteen})). But we know that ${\mathcal
L}_{\Lambda_c}\Upsilon_{e,\ov e}=0$, so that in particular
$\big({\mathcal L}_{h{\sf s}} v(c^h_e,c^{\ov h}_{\ov
e})\big)_v(0)=0$. It follows that $\ov\nabla_{\sf s}\ov\sigma =
[h{\sf s},v\ov\sigma]_v$.

For the second part, we should specify in the first place what is
meant by $v(\sigma)$: any $\sigma\in Sec(\pi)$ can be thought of
as a section of $\tpi$ and then
$v(\sigma)(e)=v(e,\vartheta_e(\sigma(\pi(e))))$. Making use of the
canonical section $\cans$ of $\pi^*\tpi$, we can write in fact
that $v(\sigma)=v(\sigma-\cans)$, where $\sigma-\cans\in
Sec(\pi^*\ov\pi)$. It is clear that $[h{\sf s},v\sigma]$ is
vertical again, and we find the properties: $\forall f \in \cinfty
M$,
\begin{eqnarray*}
[h(f{\sf s}),v\sigma]_v &=& f[h{\sf s},v\sigma]_v, \\[0mm] [h{\sf
s},v(\sigma+f\ov\sigma)]_v &=& [h{\sf
s},v(\sigma-\cans)+f v(\ov\sigma)]_v = [h{\sf s},v\sigma]_v+
f[h{\sf s},v\ov\sigma]_v + \varrho({\sf s})(f)\ov\sigma \\
&=& [h{\sf s},v\sigma]_v + f \ov\nabla_{\sf s}\ov\sigma +
\varrho({\sf s})(f)\ov\sigma.
\end{eqnarray*}
Again, these are the characterising properties of the covariant
derivative $\nabla_{\sf s}\sigma$. It follows that the operator
$L({\sf s},\sigma)=\ov\nabla_{\sf s}\sigma - [h{\sf s},v\sigma]_v$
is linear in ${\sf s}$ and affine in $\sigma$. This is the
analogue, when there are affine components involved, of an
operator $L$ being tensorial. The rest of the reasoning follows
the same pattern as before. This time, starting from an arbitrary
$e\in E_m$ and an integral curve $c^h_e$ of $h({\sf s})$ through
$e$, we put $\ov e=\sigma(\pi(e))-e$ and choose the curves
$\ov\eta=c^{\ov h}_{\ov e}$ in $\ov E$ and $\eta=c^h_e+\ov\eta$ in
$E$ to obtain a section of $\pi$ along $c^h_e$ which has
$\sigma(\pi(e))$ as initial value.\qed

{\sc Remark}. A more direct, but perhaps geometrically less
appealing proof, consists in verifying the statements of
Proposition~\ref{prop5} by a coordinate calculation. We have, for
${\sf s}={\sf s}^a(x){\sf e}_a$ and
$\sigma=e_0+\sigma^\alpha(x){\ov e}_\alpha$,
\[
\nabla_{\sf s}\sigma = [h{\sf s},v\sigma]_v = \big( \rho^i_a
\fpd{\sigma^\alpha}{x^i} + \Gamma^\alpha_{a0}(x) +
\Gamma^\alpha_{a\beta}(x)\sigma^\beta(x)\big){\sf s}^a(x) {\ov
e}_\alpha,
\]
and similarly, for $\ov\sigma={\sigma}^\alpha {\ov e}_\alpha$,
\[
\ov\nabla_{\sf s}\ov\sigma = [h{\sf s},v\ov\sigma]_v = \big(
\rho^i_a \fpd{{\sigma}^\alpha}{x^i} +
\Gamma^\alpha_{a\beta}(x){\sigma}^\beta(x)\big){\sf s}^a(x) {\ov
e}_\alpha.
\]

\begin{cor}
If $(e,{\sf v})\in \pi^*{\sf V}$ and ${\sf s}\in Sec(\tau)$,
$\sigma\in Sec(\pi)$ are sections passing through ${\sf v}$ and
$e$ respectively, then we have the following relation
\begin{equation} \label{hish}
h(e,{\sf v})= T\sigma(\varrho({\sf v})) - [h{\sf s},v\sigma](e).
\end{equation}
\end{cor}
\proof\, It was shown in \cite{MeSaMa} that a pair of operators
having the properties of covariant derivatives
$(\nabla,\ov\nabla)$ uniquely define an affine
$\varrho$-connection on $\pi$. The brackets $[h{\sf s},v\sigma]_v$
and $[h{\sf s}, v\ov\sigma]_v$ constitute such a pair (as shown
above) and according to \cite{MeSaMa}, the right-hand side of
(\ref{hish}) would then define the associated affine
$\varrho$-connection. A priori, however, there is no reason why
this would be the $h$ we started from. But the proof of
Proposition \ref{prop5} precisely guarantees now that it must be
the $h$ we started from, and hence we have (\ref{hish}). \qed

There is of course a similar formula for $\ov h$, which reads,
\begin{equation}
 \label{ovhisovh}
\ov h(\ov e,{\sf v})= T\ov\sigma(\varrho({\sf v})) - [h{\sf
s},v\ov\sigma](\ov e).
\end{equation}

As an immediate benefit of the formulas (\ref{mm}), we can obtain
an explicit defining relation now for the extension of the
operators $(\nabla,\ov\nabla)$ to a covariant derivative
$\tilde\nabla$ on $Sec(\tpi)$. Each $\tilde\sigma\in Sec(\tpi)$ is
either of the form $f\sigma$, with $f\in\cinfty M$ and $\sigma\in
Sec(\pi)$, or of the form $\ov\sigma$, for some $\ov\sigma\in
Sec(\ov\pi)$. Then, $\tilde\nabla_{\sf s}\tilde\sigma$ is defined
in \cite{MeSaMa} by one or the other of the following relations
\[
\tilde\nabla_{\sf s}\tilde\sigma=f\nabla_{\sf s}\sigma +
\varrho({\sf s})(f)\sigma, \qquad \mbox{or} \qquad
\tilde\nabla_{\sf s}\tilde\sigma = \ov\nabla_{\sf s}\ov\sigma.
\]
\begin{cor}
A unifying formula for the computation of $\tilde\nabla_{\sf
s}\tilde\sigma$ is given by
\begin{equation} \label{decnablatilde}
\tilde\nabla_{\sf s}\tilde\sigma = [h{\sf s},v\tilde\sigma]_v +
\varrho({\sf s})(\langle\tilde\sigma,e^0\rangle)\cans.
\end{equation}
\end{cor}
\proof\, In the first case we have $\tilde\nabla_{\sf
s}\tilde\sigma= f[h{\sf s},v\sigma]_v + \varrho({\sf
s})(f)\sigma$. Using the properties that $\sigma$, regarded as
section of $\pi^*\tpi$, can be written as $\sigma =
\vartheta(\sigma)+\cans$, and that
$\big(v(\sigma)\big)_v=\vartheta(\sigma)$, this expression can be
rewritten as $\tilde\nabla_{\sf s}\tilde\sigma=[h{\sf
s},v(f\sigma)]_v + \varrho({\sf s})(f)\cans$, which is of the form
(\ref{decnablatilde}) since $\langle\tilde\sigma,e^0\rangle = f$
in this case. In the second case, we have $\tilde\nabla_{\sf
s}\tilde \sigma = [h{\sf s},v{\ov\sigma}]_v$, which is immediately
of the right form since $\langle\tilde\sigma,e^0\rangle=0$ now.
\qed

The representation (\ref{decnablatilde}) of $\tilde\nabla_{\sf
s}\tilde\sigma$ is exactly the decomposition (\ref{decomp}) of
$\tilde\nabla_{\sf s}\tilde\sigma$, regarded as section of
$\pi^*\tpi$. One should not forget, of course, that such a
decomposition somehow conceals part of the information in case the
section under consideration, as is the case with
$\tilde\nabla_{\sf s}\tilde\sigma$ here, is basic, in the sense
that it is actually a section of $\tpi:\tilde E \rightarrow M$.

\section{Generalised connections of Berwald type}

Within the framework of the classical theory of connections on a
tangent bundle $\tau_M:TM \rightarrow M$ (or more generally a
vector bundle over $M$), it is well known that an arbitrary
horizontal distribution or non-linear connection has a kind of
linearisation \cite{Vilms}. For the case of the tangent bundle,
for example, this induced linear connection can be interpreted in
different equivalent ways: as a linear connection on
$T(TM)\rightarrow TM$ (see for example \cite{Szil}), or as a
linear connection on the vertical bundle $V(TM)\rightarrow TM$
(see e.g.\ \cite{Bejancu}), or perhaps most efficiently as a
connection on the pullback bundle $\tau_M^*TM\rightarrow TM$. An
interesting geometrical characterisation of this so-called
Berwald-type connection, in its pullback bundle version, was given
by Crampin in \cite{Crampin}. Our generalisation to a
time-dependent set-up on jet bundles \cite{MeSa} revealed that
there is a certain liberty in fixing the time-component of the
connection, though two particular choices come forward in a rather
natural way via a direct defining relation of the covariant
derivative. This kind of gauge freedom in fixing the connection
has everything to do with the affine nature of the first-jet
bundle. We shall now explore to what extent a $\varrho$-connection
on an affine bundle, in the general picture of Section~2, has a
kind of induced linearisation, and we intend to unravel in that
process the origin of the two specific choices for fixing the
connection, as described in \cite{MeSa}.

The case of our time-dependent model in \cite{MeSa} fits within
the present scheme as follows: ${\sf V}=TM$, $\varrho$ is the
identity and $E$ is the first-jet bundle of $M\rightarrow \R$. The
induced linear connection then is a connection on the pullback
bundle $\pi^*\tpi$, i.e.\ a covariant derivative operator
$\nabla_\xi X$, where $\xi$ is a vector field on $E$ and $X$ a
section of $\pi^*\tpi$. To define $\nabla_\xi X$, it suffices to
specify separately the action of horizontal and vertical vector
fields, where ``horizontality'' is defined of course via the
non-linear connection one starts from. In the more general
situation of a $\varrho$-connection, however, horizontality of
vector fields on $E$ is not an unambiguous notion, in the sense
that ${\rm{Im}}\, h$ may not provide a full complement of the set
of vertical vectors and may even have a non-empty intersection
with this set (see \cite{FransBavo}). As said in
Section~\ref{basicsetupsection}, we do have a direct complement
for the vertical sections of $\pi^1:\prolr E \rightarrow E$. So
the right way to look here for a linear connection on $\pi^*\tpi$
is as a $\varrho^1$-connection.

The linear $\varrho^1$-connection on $\pi^*\tpi$ will actually be
generated by an affine $\varrho^1$-connection on $\pi^*E
\rightarrow E$. However, as long as we let ${\sf V}\rightarrow M$
be any vector bundle, not related to $E\rightarrow M$ and without
the additional structure of a Lie algebroid, there is no bracket
of sections of $\tau$ or $\pi^1$ available. We should, therefore,
not expect to discover easily direct defining relations. Instead,
we shall approach the problem of detecting corresponding
$\varrho^1$-connections on $\pi^*E$ via their covariant derivative
operators, for which we will use the results of Proposition
\ref{prop5} as one of the sources of inspiration. Following the
lead of Crampin's approach in \cite{Crampin}, the other source of
inspiration should come from understanding the details of possible
rules of parallel transport, into which subject we will enter now
first.

Recall that the concept of parallel transport in $E$, i.e. the
construction of the horizontal lift $c^h_e$ of a
$\varrho$-admissible curve in ${\sf V}$, exists for any
$\varrho$-connection $h$. If $h$ is affine, we know that for
horizontal lifts which start at $e$ and $e_1=e+\ov e$ at an
initial time $a$, we have at any later time $b$ that $c^h_{e+\ov
e}(b)=c^h_e(b)+c^{\ov h}_{\ov e}(b)$. In addition, $v(e,\ov e)$
identifies the couple $(e,\ov e)\in\pi^*\ov E$ with a vertical
tangent vector to $E$, and we have seen that the evolution to the
vector $v(c^h_{e},c^{\ov h}_{\ov e})$ is just Lie transport along
$c^h_e$. If $h$ is {\em not} affine, Lie transport of a vertical
vector along $c^h_e$ still exists and one could somehow reverse
the order of thinking to use that for defining an affine action on
fibres of $E$. To be specific, writing $Y_{e,\ov e}(a)=v(e,\ov e)$
for the initial vertical vector and considering its Lie transport,
defined as before by $Y_{e,\ov e}(u)=T\phi^h_{u,a}\big(Y_{e,\ov
e}(a)\big)$, we get the following related actions on $\pi^*\ov E$
and $\pi^*E$:
\begin{equation}
\begin{array}{l}
(e,\ov e) \mapsto (c^h_e,p_{\ov E}\big((Y_{e,\ov e})_v\big)), \\[2mm]
(e,e_1) \mapsto (c^h_e,c^h_e + p_{\ov
E}\big((Y_{e,e_1-e})_v\big)).
\end{array}
\end{equation}
We will refer to this as the {\em affine action on $\pi^*E$ by Lie
transport along horizontal curves}. Obviously, when $h$ is not
affine, the image of $(e,e_1)$ under this affine action will not
be $(c^h_e,c^h_{e_1})$.

The question which arises now is whether there are natural ways
also to define an affine action on $\pi^*E$ along vertical curves,
i.e.\ curves in a fixed fibre $E_m$ of $E$. Let $c^v_e$ denote an
arbitrary curve through $e$ in the fibre $E_m$ ($m=\pi(e)$). It
projects onto the constant curve $c_m: u \mapsto c_m(u)=m$ in $M$.
A curve in ${\sf V}$ which has the same projection (and actually
is $\varrho$-admissible) can be taken to be ${\sf o}_m: u \mapsto
{\sf o}_m(u) = (m,{\sf o}_m)$. ${\dot c}^v_e$ is a curve in $TE$
which projects onto $c^v_e$ and has the property $T\pi({\dot
c}^v_e)=0$. By analogy with earlier constructions, we define a new
curve $\V{\dot c}_e$ in $\prolr E$, determined by
\begin{equation}
\V{\dot c}_e : = ({\sf o}_m,{\dot c}^v_e). \label{krulleke}
\end{equation}
Obviously, by construction, we have that $\pi^1\circ \V{\dot c}_e
= c^v_e$ and $\varrho^1\circ \V{\dot c}_e = {\dot c}^v_e$, so that
$\V{\dot c}_e$ is $\varrho^1$-admissible.

Let us now address the problem of defining a transport rule in
$\pi^*E$ along curves $c^v_e$. Remember that for the horizontal
curves, we described such a transport rule by looking first at the
way vertical tangent vectors can be transported. For the transport
of vertical tangent vectors within a fixed fibre, the usual
procedure is to take simple translation (this is sometimes called
{\em complete parallelism}). Thus, starting from a point
$(e,e_1)\in \pi^*E$, to which we want to associate first a
vertical tangent vector, we think of $(e,e_1)$ as belonging to
$\pi^*\tilde E$ and consider
$v(e,e_1)=v(e,\vartheta_e(e_1))=v(e,e_1-e)$. Its parallel
translate along a curve $c^v_e$ is $v(c^v_e,e_1-e)$ which can be
identified with $(c^v_e,e_1-e)\in\pi^*\ov E$. But it makes sense
to associate with this a new element of $\pi^*E$ as well, in
exactly the same way as we did it for horizontal curves. We thus
arrive at the following action on $\pi^*\ov E$ and $\pi^*E$
\begin{equation}
\begin{array}{l}
(e,\ov e) \mapsto (c^v_e,\ov e), \\[2mm]
(e,e_1) \mapsto (c^v_e,c^v_e +e_1-e).
\end{array} \label{fwactionovE}
\end{equation}
It could be described as a {\em vertical affine action by
translation in $\pi^*\ov E$}.

There is, however, another way of transporting points in $\pi^*E$
along a curve of type $c^v_e$, which is in fact the most obvious
one if one does insist on having a link with a transport rule of
vertical tangent vectors via the vertical lift operator. It is
obtained by looking at the action
\begin{equation}
\begin{array}{l}
(e,\ov e) \mapsto (c^v_e,\ov e), \\[2mm]
(e,e_1) \mapsto (c^v_e,e_1),
\end{array} \label{fwactionE}
\end{equation}
and could be termed as a {\em vertical affine action by
translation in $\pi^*E$}.

Given an arbitrary $\varrho$-connection $h$ on the affine bundle
$\pi: E\rightarrow M$, we now want to construct an  induced
$\varrho^1$-connection $h^1$ on the affine bundle $\pi^*\pi:
\pi^*E\rightarrow E$ through the identification of suitable
covariant derivative operators $\cd{}$ and $\cdov{}$. That is to
say, we should give a meaning to things like $\cd{\mathcal Z}X$
and $\cdov{\mathcal Z}\ov X$, for ${\mathcal Z}\in Sec(\pi^1)$,
$X\in Sec(\pi^*\pi)$, $\ov X\in Sec(\pi^*\ov\pi)$. As explained in
Section~2, every ${\mathcal Z}$ has a unique decomposition in the
form ${\mathcal Z}=\H{\sf X}+\V{\ov Y}$, with ${\sf X}\in
Sec(\pi^*\tau)$, $\ov Y \in Sec(\pi^*\ov\pi)$. These in turn are
finitely generated (over $\cinfty E$) by {\em basic sections},
i.e. sections of $\tau$ and of $\ov\pi$, respectively. The same is
true for the sections $X$ or $\ov X$ on which $\cd{\mathcal Z}$
and $\cdov{\mathcal Z}$ operate. This means that, for starting the
construction of covariant derivatives, we must think of a defining
relation for $\cd{\H{\sf s}}\sigma$, $\cd{\V{\ov \eta}}\sigma$,
$\cdov{\H{\sf s}}\ov\sigma$, $\cdov{\V{\ov \eta}}\ov\sigma$, with
${\sf s}\in Sec(\tau)$, $\sigma\in Sec(\pi)$,
$\ov\eta,\ov\sigma\in Sec(\ov\pi)$. The expectation is, since we
look for a $\cd{}$ and $\cdov{}$, that $h^1$, as a kind of
linearisation of $h$, will be an affine connection and so, in the
particular case that the given $h$ is affine, it should
essentially reproduce a copy of itself. Therefore, the first idea
which presents itself is to set
\begin{equation}
\cd{\H{\sf s}}\sigma = [h{\sf s},v\sigma]_v, \quad \cdov{\H{\sf
s}}\ov\sigma = [h{\sf s},v\ov\sigma]_v, \quad
\cd{\V{\ov\eta}}\sigma=\cdov{\V{\ov\eta}}\ov\sigma=0.
\label{local1}
\end{equation}
The first point in the proof of Proposition \ref{prop5} did not
rely on the assumption of $h$ being affine, so we know that these
formulas at least are consistent with respect to the module
structure over $\cinfty M$. We then extend the range of the
operators $\cd{}$ and $\cdov{}$ in the obvious way, by the
following three rules: for every $F\in \cinfty E$, we put
\begin{eqnarray}
&& \cd{F\H{\sf s}}\sigma =  F \cd{\H{\sf s}}\sigma = F [h{\sf
s},v\sigma]_v, \quad \cdov{F\H{\sf s}}\ov\sigma = F \cdov{\H{\sf
s}}\ov\sigma = F [h{\sf s},v\ov\sigma]_v,
\label{local2} \\
&& \cd{F\V{\ov\eta}}\sigma=\cdov{F\V{\ov\eta}}\ov\sigma=0,
\label{local3}
\end{eqnarray}
which suffices to know what $\cd{\mathcal Z}\sigma$ and
$\cdov{\mathcal Z}\ov \sigma$ mean for arbitrary ${\mathcal Z}\in
Sec(\pi^1)$, and finally we put
\begin{eqnarray}
\cdov{\mathcal Z}(F\ov\sigma) = F\cdov{\mathcal Z}\ov\sigma
+ \varrho^1({\mathcal Z})(F)\ov\sigma, \label{local4} \\
\cd{\mathcal Z}(\sigma+F\ov\sigma) = \cd{\mathcal Z}\sigma+
F\cdov{\mathcal Z}\ov\sigma + \varrho^1({\mathcal Z})(F)\ov\sigma,
\label{local5}
\end{eqnarray}
which suffices to give a meaning to all $\cd{\mathcal Z}X$ and
$\cdov{\mathcal Z}\ov X$. Our operators satisfy by construction
all the necessary requirements for defining an affine
$\varrho^1$-connection $h^1$.

It is worthwhile to observe that for the covariant derivatives of
general $X\in Sec(\pi^*\pi)$ and $\ov X\in Sec(\pi^*\ov\pi)$, we
still have an explicit formula at our disposal when ${\mathcal Z}$
is of the form $\H{\sf s}$, with ${\sf s}$ basic. This follows
from the fact that $\varrho^1(\H{\sf s})=h({\sf s})$, so that
\begin{eqnarray}
\cd{\H{\sf s}}(\sigma+F\ov\sigma)&=& [h{\sf s},v\sigma]_v+F[h{\sf
s},v\ov\sigma]_v +
\varrho^1(\H{\sf s})(F)\ov\sigma, \nonumber \\
&=& [h{\sf s},v(\sigma+F\ov\sigma)]_v, \label{local6}
\end{eqnarray}
and likewise for $\cd{\H{\sf s}}\ov X$.

The next point on our agenda is to understand what parallel
transport means for the affine connection $(\cd{},\cdov{})$, or
even better, to show that it is uniquely characterised by certain
features of its parallel transport. The general idea of parallel
transport is clear, of course: starting from any
$\varrho^1$-admissible curve $c^1$ in $\prolr E$, its horizontal
lift is a curve $\psi^1$ in $\pi^*E$ having the same projection
$\psi^1_{\scriptscriptstyle E}=c^1_{\scriptscriptstyle E}$ in $E$
and satisfying $\cd{c^1}\psi^1=0$; image points of $\psi^1$ then
give parallel translation by definition. Now, $\psi^1$ is
essentially a pair of curves in $E$ having the same projection in
$M$, so the determination of $\psi^1$ is a matter of constructing
a second curve in $E$ having the same projection in $M$ as
$c^1_{\scriptscriptstyle E}$. It will be sufficient to focus on
$\varrho^1$-admissible curves of the form $\H{\dot c}_e$ and
$\V{\dot c}_e$, for which the corresponding projections on $E$ are
curves of the form $c^h_e$ and $c^v_e$, respectively, and to
consider curves $\psi^1$ which come from the restriction of
sections of $\pi^*\pi$ to $c^h_e$ or $c^v_e$. To simplify matters
even further, we can use basic sections ${\sf s}\in Sec(\tau)$ to
generate horizontal curves, because we know from Lemma
\ref{Lemma1} that the integral curves of $h({\sf s})\in
\vectorfields E$ are horizontal lifts. Vertical curves, of course,
can be generated as integral curves of vertical vector fields.

\begin{prop}
Let ${\sf s}\in Sec(\tau)$, $\ov Y \in Sec(\pi^*\ov\pi)$ be
arbitrary. Denote the integral curves of $h({\sf s})$ and $v\ov Y$
through a point $e$ by $c^h_e$ and $c^v_e$ and consider their
lifts to $\varrho^1$-admissible curves $\H{\dot c}_e$ and $\V{\dot
c}_e$ in $\prolr E$. $(\cd{},\cdov{})$ is the unique affine
$\varrho^1$-connection on $\pi^*\pi$ with the properties
\begin{itemize}
\item[(i)] Parallel transport along $\H{\dot c}_e$ is the affine
action on $\pi^*E$ by Lie transport along horizontal curves.
\item[(ii)] Parallel transport along $\V{\dot c}_e$ is the
vertical affine action by translation in  $\pi^*E$.
\end{itemize} \label{prop6}
\end{prop}
\proof\, Recall that $\H{\sf s}\in Sec(\pi^1)$ is defined at each
$e\in E$ by $\H{\sf s}(e)=\big({\sf s}(\pi(e)),h({\sf
s})(e)\big)$, so that at each point along an integral curve
$c^h_e$ of $h({\sf s})$, we have
\[
\H{\sf s}(c^h_e(u))=\big({\sf s}\circ\pi\circ c^h_e(u), {\dot
c}^h_e(u)\big)=\H{\dot c}_e(u).
\]
Let now $X$ be an arbitrary section of $\pi^*\pi$ and put
$\psi^1(u)=X(c^h_e(u))$, which defines a curve in $\pi^*E$
projecting onto $c^h_e$ in $E$. We have
\[
\big(\cd{\H{\dot c}_e}\psi^1\big)(u)=\cd{\H{\dot c}_e(u)}X =
\cd{\H{\sf s}(c^h_e(u))}X = \big(\cd{\H{\sf
s}}X\big)\big(c^h_e(u)\big).
\]
If such curve is required to govern parallel transport in
$\pi^*E$, we must have $\big(\cd{\H{\dot c}_e}\psi^1\big)(u)=0$,
$\forall u$. This implies that $\forall {\sf s}\in Sec(\tau)$,
$\forall X \in Sec(\pi^*\pi)$, $\cd{\H{\sf s}}X$ should be zero
along integral curves of $h({\sf s})\in\vectorfields E$. By the
remark about the explicit formula (\ref{local6}) for $\cd{\H{\sf
s}}X$ and with $v(X)(c^h_e(u))=v(\psi^1(u))$, which defines a
vertical vector field along $c^h_e$, this requirement is further
equivalent to ${\mathcal L}_{h({\sf s})}v(\psi^1)=0$, which is
precisely the characterisation of Lie transport. The same
arguments apply to $\cdov{}$ and show that our $(\cd{},\cdov{})$
has the property (i).

With $\ov Y\in Sec(\pi^*\ov\pi)$, $\V{\ov Y}\in Sec(\pi^1)$ is
such that $\varrho^1(\V{\ov Y})$ is a vertical vector field on
$E$. Hence, its integral curves are curves of the form $c^v_e$ in
a fixed fibre $E_{\pi(e)}$ and we have from (\ref{krulleke}):
\[
\V{\ov Y}\big(c^v_e(u)\big)=\big({\sf o}_{\pi(e)},{\dot
c}^v_e(u)\big) = \V{\dot c}_e(u).
\]
Let $X$ again be an arbitrary section of $\pi^*\pi$ and put this
time $\psi^1(u)=X(c^v_e(u))$, which defines a curve in $E\times_M
E$ projecting onto $c^v_e$ for its first component. We wish to
show that if $\psi^1(u)$ rules parallel transport along $\V{\dot
c}_e$, it is necessarily a curve which is constant in its second
component. We have
\[
\big(\cd{\V{\dot c}_e}\psi^1\big)(u)=\cd{\V{\dot c}_e(u)}X =
\cd{\V{\ov Y}(c^v_e(u))}X = \big(\cd{\V{\ov
Y}}X\big)\big(c^v_e(u)\big).
\]
There is no explicit formula available for $\cd{\V{\ov Y}}X$.
However, $X$ is locally of the form $X=\sigma+F_i{\ov\sigma}_i$,
with $\sigma$, ${\ov\sigma}_i$ basic sections and $F_i\in\cinfty
E$. It then follows that $\cd{\V{\ov Y}}X=\varrho^1(\V{\ov
Y})(F_i){\ov\sigma}_i$ and the requirement $\cd{\V{\ov
Y}}X(c^v_e(u))=0$ implies that the $F_i$ must be the first
integrals of $\varrho^1(\V{\ov Y})$. In turn this means that the
value $X(c^v_e(u))$ is constant. This way we see that the affine
connection $(\cd{},\cdov{})$ also has property (ii).

That properties (i) and (ii) uniquely fix the connection is easy
to see, because the above arguments show that they impose in
particular that $\cd{\H{\sf s}}\sigma=[h{\sf s},v\sigma]_v$ and
$\cd{\V{\ov\eta}}\sigma=0$ (and similarly for $\cdov{}$), for
basic ${\sf s},\sigma$ and $\ov\eta$. And these are exactly the
defining relations (\ref{local1}) from which our couple
$(\cd{},\cdov{})$ was constructed.\qed

We have seen earlier on that there is a second interesting
transport rule along vertical curves and would like to discover
now what modifications to the affine connection must be made to
have this other rule as vertical parallel transport. We are
referring here to the action (\ref{fwactionovE}) for which the
curve starting at some $(e,e_1)\in E_m\times E_m$ is of the form
\[
u\overset{\psi^1}{\mapsto} \big(c^v_e(u),e_1+c^v_e(u)-e\big) =
\big(c^v_e(u),c^v_e(u)+e_1-e\big).
\]
Now, it is easy to identify a section of $\pi^*\pi$ which along
$c^v_e$ coincides with this curve. Indeed, choosing a basic
section $\ov\sigma\in Sec(\ov\pi)$ which at $m=\pi(e)$ coincides
with $e_1-e$, we are simply looking at the restriction to $c^v_e$
of $\cans + \ov\sigma$, where $\cans$ here denotes the identity
map on $E$.

Let $(\cdhat{}, \cdhatov{})$ denote the affine connection we are
looking for now and which clearly will coincide with
$(\cd{},\cdov{})$ for its ``horizontal action''. If, as before,
$\ov Y\in Sec(\pi^*\ov\pi)$ generates the vertical vector field
$v{\ov Y}$ whose integral curves are the $c^v_e$, the above
$\psi^1$ will produce parallel transport, provided we have
\[
\big(\cdhat{\V{\dot c}_e}\psi^1\big)(u)= \cdhat{\V{\ov
Y}}(\cans+\ov\sigma)(c^v_e(u))=0.
\]
Since this must hold for each $\V{\ov Y}$ and, for every fixed
$\V{\ov Y}$ also for all $\ov\sigma$, this is equivalent to
requiring that $\cdhat{\V{\ov Y}}\cans=0$ and $\cdhatov{\V{\ov
Y}}\ov\sigma=0$, $\forall \ov Y,\ov\sigma$. In fact, in view of
the linearity in $\ov Y$, we actually obtain the conditions
\begin{equation}
\cdhat{\V{\ov\eta}}\cans=0 \qquad \mbox{and} \qquad
\cdhatov{\V{\ov\eta}}\ov\sigma=0, \label{local7}
\end{equation}
for all basic $\ov\sigma$ and $\ov\eta$. It is interesting to
characterise this completely by properties on basic sections,
because the extension to a full affine connection on $\pi^*\pi$
then follows automatically. If $\sigma$ is an arbitrary basic
section of $\pi^*\pi$, it can be decomposed (see(\ref{decomp})) in
the form $\sigma=\cans+ \vartheta(\sigma)$, whereby
$\vartheta(\sigma)(e)= \big(e,\sigma(\pi(e))-e\big)$. Clearly,
$\sigma(\pi(e))-e$, as an element of $\ov E$, has components which
are linear functions of the fibre coordinates of $e$, in such a
way that when acted upon by the vector field
$\varrho^1(\V{\ov\eta})$, we will obtain $-{\ov\eta}$. It follows
that $\cdhatov{\V{\ov\eta}}\vartheta(\sigma) =-\ov\eta$ and
therefore that
\begin{equation}
\cdhat{\V{\ov\eta}}\cans =0 \qquad \Longleftrightarrow \qquad
\cdhat{\V{\ov\eta}}\sigma=-\ov\eta, \,\, \forall \sigma \in
Sec(\pi). \label{local8}
\end{equation}
This way, we have detected an alternative way for defining an
affine $\varrho^1$-connection $(\cdhat{}, \cdhatov{})$ on
$\pi^*\pi$. Compared to (\ref{local1}), its defining relations are
\begin{equation}
\cdhat{\H{\sf s}}\sigma = [h{\sf s},v\sigma]_v, \quad
\cdhatov{\H{\sf s}}\ov\sigma = [h{\sf s},v\ov\sigma]_v, \quad
\cdhat{\V{\ov\eta}}\sigma=-\ov\eta,\quad
\cdhatov{\V{\ov\eta}}\ov\sigma=0. \label{local9}
\end{equation}
We can further immediately draw the following conclusion about its
characterisation
\begin{prop}
With the same premises as in Proposition \ref{prop6},
$(\cdhat{},\cdhatov{})$ is the unique affine
$\varrho^1$-connection on $\pi^*\pi$ with the properties
\begin{itemize}
\item[(i)] Parallel transport along $\H{\dot c}_e$ is the affine
action on $\pi^*E$ by Lie transport along horizontal curves.
\item[(ii)] Parallel transport along $\V{\dot c}_e$ is the
vertical affine action by translation in  $\pi^*\ov E$.
\end{itemize}\label{prop7}
\end{prop}

We will refer to the connections $(\cd{},\cdov{})$ and
$(\cdhat{},\cdhatov{})$, as well as their extensions $\cdtilde{}$
and $\cdhattilde{}$ as {\em Berwald-type connections}. For
completeness, we list their defining relations here in
coordinates.
\[
\begin{array}{ll}  \displaystyle{\cd{{\mathcal H}_a} e_0 =
\left(\Gamma^\gamma_a-y^\beta\fpd{\Gamma^\gamma_a}{y^\beta}\right)
{\ov e}_\gamma}, & \qquad \qquad \displaystyle{\cd{{\mathcal V}_\alpha} e_0 = 0}, \\[5mm]
\displaystyle{\cdov{{\mathcal H}_a} {\ov e}_\beta =
\fpd{\Gamma^\gamma_a}{y^\beta} {\ov e}_\gamma}, & \qquad \qquad
\displaystyle{\cdov{{\mathcal V}_\alpha} {\ov e}_\beta = 0}
\end{array}
\]
and
\[
\begin{array}{ll}  \displaystyle{\cdhat{{\mathcal H}_a} e_0 =
\left(\Gamma^\gamma_a-y^\beta\fpd{\Gamma^\gamma_a}{y^\beta}\right)
{\ov e}_\gamma}, & \qquad \qquad \displaystyle{\cdhat{{\mathcal V}_\alpha}
e_0 = -{\ov e}_\alpha}, \\[5mm]
\displaystyle{\cdhatov{{\mathcal H}_a} {\ov e}_\beta =
\fpd{\Gamma^\gamma_a}{y^\beta} {\ov e}_\gamma}, & \qquad \qquad
\displaystyle{\cdhatov{{\mathcal V}_\alpha} {\ov e}_\beta = 0}.
\end{array}
\]

\section{The case of affine Lie algebroids and the canonical connection
associated to a pseudo-\sode}

It is now time to relate the new and quite general results of the
preceding sections to some of our earlier work. Recall that our interest
in affine bundles comes in the first place from the geometrical study
of time-dependent second-order equations and the analysis of
Berwald-type connections in that context \cite{MeSa}. Secondly, once
the potential relevance for applications of Lagrangian systems on
Lie algebroids became apparent, we were led to explore a time-dependent
generalisation of such systems and thus arrived at the introduction
and study of affine Lie algebroids \cite{SaMeMa,MaMeSa}. Notice that
Lagrangian systems on algebroids are particular cases of pseudo-\sode s,
but the concept of a pseudo-\sode\ in itself, strictly speaking, does
not require the full structure of a Lie algebroid.

So, let us start by looking at pseudo-\sode s on the affine bundle $E$,
which are essentially vector fields with the property that all the
integral curves are $\rho$-admissible. In saying that, we are in fact
assuming that the anchor map has $E$ in its domain.
In this section, therefore, the starting point is that we have an
affine bundle map $\rho:E\rightarrow TM$ at our disposal with associated
linear map $\ov \rho: \ov E \rightarrow TM$.
These maps can be extended to a
linear map $\tilde\rho: \tilde E \rightarrow TM$ as follows:
for every $\tilde e\in\tilde{E}_m$,
making a choice of an element $e\in E_m$, we have a representation
of the form $\tilde e=\lambda e +\ov e$ and define
$\tilde\rho(\tilde e)$ by
\[
\tilde\rho(\tilde e)= \lambda \rho(e) + {\ov \rho}({\ov e}).
\]
One easily verifies that this construction does not depend on the
choice of $e$. In coordinates:
\[
\rho:(x^i,y^\alpha)\mapsto
\left(\rho^i_\alpha(x)y^\alpha + \rho^i_0(x)\right)\fpd{}{x^i},
\qquad \trho: (x^i,y^0,y^\alpha)\mapsto \left(\rho^i_\alpha
y^\alpha + \rho^i_0y^0\right)\fpd{}{x^i}.
\]
In what follows $\trho$ plays the role of the anchor map
$\varrho$ we had before. This means in particular that the vector
bundle $\tau: {\sf V} \rightarrow M$ from now on is taken to be
the bundle $\tpi: \tilde E\rightarrow M$. Now, pseudo-\sode s can
be regarded also as sections of the prolonged bundle
$\pi^1: \proltr E\rightarrow E$ (rather than as vector fields on $E$).
The difference in interpretation is easy to understand from the
basic constructions explained in Section~2. Indeed, we have seen there
that there is a natural vertical lift operator
$v: Sec(\pi^*\tpi)\rightarrow \vectorfields E$, which extends to an
operator $\V{}: Sec(\pi^*\tpi)\rightarrow Sec(\pi^1)$ via (\ref{vlift}).
Combining this vertical lift with the projection $j:\proltr E
\rightarrow \pi^*\tilde E$, gives rise to the map $S=\V{}\circ j$,
called the {\em vertical endomorphism on $Sec(\pi^1)$}.
For ${\mathcal Z}={\zeta}^0\baseX
0 + \zeta^\alpha\baseX \alpha + Z^\alpha \baseV\alpha$, we have
\[
S({\mathcal Z}) = (\zeta^\alpha-\zeta^0y^\alpha)\baseV\alpha.
\]
An elegant definition of the concept of pseudo-\sode\  then goes
as follows.

\begin{dfn}
A pseudo-\sode\ is a section $\Gamma$ of $\pi^1$ such that
$S(\Gamma)=0$ and $\langle\Gamma,{\mathcal X}^0\rangle = 1$.
\end{dfn}
In coordinates, $\Gamma$ is of the form
\begin{equation}
\Gamma=\baseX 0 + y^\alpha \baseX \alpha + f^\alpha \baseV \alpha.
\label{pseudocoord}
\end{equation}
It is not immediately clear whether a pseudo-\sode\ comes with
a canonically associated (non-linear) $\trho$-connection in this
general setting. However, as mentioned already in \cite{Me}, the
construction of a connection becomes quite obvious when we have
the additional structure of a Lie algebroid.

So, assume now we have an affine Lie algebroid structure on $\pi$,
which can most conveniently be seen as a Lie algebroid on $\tpi$ with anchor
$\trho$, which is such that the bracket of two sections
of $\pi$ (regarded as sections of $\tpi$) is a
section of $\ov\pi$ (also considered as a section of
$\tpi$). In coordinates, there exist
structure functions $C^\gamma_{\alpha\beta}$ and
$C^\gamma_{0\beta}$ on $M$ such that
\[
[e_\alpha,e_\beta] = C^\gamma_{\alpha\beta}(x) e_\gamma
\qquad\mbox{and}\qquad [e_0,e_\beta] =
C^\gamma_{0\beta}(x)e_\gamma.
\]
We have shown in \cite{MaMeSa} that such a Lie algebroid on $\tpi$ can
be prolonged to a Lie algebroid on
$\tpi^1$ with anchor $\trho^1$. In coordinates:
\begin{equation}
\begin{array}{l}
[\baseX \alpha, \baseX \beta] = C^\gamma_{\alpha\beta}\baseX
\gamma, \qquad [\baseX 0, \baseX \beta] = C^\gamma_{0\beta}\baseX
\gamma, \qquad [\baseV \alpha, \baseX \beta] = 0,\\[2mm]
[\baseX 0, \baseV \beta] = 0, \qquad\qquad[\baseV \alpha, \baseV
\beta] = 0.
\end{array}\label{extalg}
\end{equation}

The Lie algebroid structure provides us with an exterior derivative; we
use the standard notation $d_\Gamma$ for the commutator $[i_\Gamma,d]$,
which plays the role of Lie derivative and extends, as a degree zero
derivation, to tensor fields of any type.

Now, one way of pinning down a $\trho$-connection on $\pi$
consists in identifying its horizontal projector $\subH P$ (and then
$\subV P = I - \subH P$).

\begin{prop}
If $\Gamma$ is a pseudo-\sode\, on an affine Lie algebroid $\pi$,
then the operator
\begin{equation}
\subH P = \frac{1}{2}\Big(I-d_\Gamma S + {\mathcal X}^0 \otimes
\Gamma \Big) \label{PHpseudo}
\end{equation}
defines a horizontal projector on $Sec(\pi^1)$ and hence a
$\trho$-connection on $\pi$.
\end{prop}
\proof\, The proof follows the lines of the classical one for
time-dependent mechanics (see \cite{CPT} or \cite{Crampinbook}).
Since it was largely omitted in \cite{Me}, we give a brief sketch
of one possibility to proceed here. For $\tilde \sigma \in Sec(\tpi)$,
define the horizontal lift $\H{\tilde\sigma}\in Sec(\pi^1)$ by
\begin{equation}
\H{\tilde\sigma}=
\frac{1}{2}\Big({\tilde\sigma}^{\scriptscriptstyle C}+ \langle
\tilde \sigma, e^0\rangle \Gamma - [\Gamma,\V{\tilde\sigma}] \Big),
\label{Hpseudo}
\end{equation}
where ${\tilde\sigma}^{\scriptscriptstyle C}$ is the complete lift,
as defined in \cite{MaMeSa}. It is easy to see that this behaves
tensorially for multiplication by basic functions and that
$\H{\tilde\sigma}$ projects onto $\sigma$. Hence, extending the
horizontal lift to $Sec(\pi^*\tpi)$ by imposing linearity for
multiplication by functions on $E$, we obtain a splitting of the
short exact sequence (\ref{shortexactgeneral}) for the present
situation. This in fact concludes the proof of the existence
of a $\trho$-connection, but it is interesting to verify further
the explicit formula for $\subH P$. One can, for example, compute
the Lie algebroid brackets $[\Gamma,\V{\ov\sigma}]$ and
$[\Gamma,\H{\ov\sigma}]$ for $\ov\sigma \in Sec(\ov\pi)$, from
which it then easily follows (using also the properties
$S(\H{\tilde\sigma})=\V{\ov\sigma}$ and $S(\V{\ov\sigma})=0$), that
$d_\Gamma S(\V{\ov\sigma})= \V{\ov\sigma}$,
$d_\Gamma S(\H{\ov\sigma})=- \H{\ov\sigma}$ and
$d_\Gamma S(\Gamma)=0$. The verification that $\subH P$ is a projection
operator and that $\subH P(\H{\tilde\sigma}) = \H{\tilde\sigma}$ then
is immediate. \qed

The connection coefficients of the pseudo-\sode\ connection are given by
\begin{eqnarray}
\Gamma^\alpha_0 &=& - f^\alpha +
\frac{1}{2}y^\beta\Big(\fpd{f^\alpha}{y^\beta}+C^\alpha_{0\beta}\Big)
=-f^\alpha-y^\beta\Gamma^\alpha_\beta \qquad \label{conncoeff1} \\
\Gamma^\alpha_\beta &=& - \frac{1}{2}\Big(\fpd{f^\alpha}{y^\beta}+
y^\gamma C^\alpha_{\gamma\beta} + C^\alpha_{0\beta} \Big).
\label{conncoeff2}
\end{eqnarray}

Briefly, the particular case of a Lagrangian system
on the affine Lie algebroid $\pi$ is obtained as follows.
Let $L$ be a function on $E$ and consider the 1-form
$\theta_L= dL \circ S + L {\mathcal X}^0$. If $\omega_L= d\theta_L$
has maximal rank at every point, i.e.\ when $L$ is said to be regular,
there exists a unique pseudo-\sode\, such that $i_\Gamma\omega_L=0$.
In that case, the functions $f^\alpha$ which determine the connection
coefficients are given by
\[
f^\alpha = g^{\alpha\beta}\left( \rho^i_\beta \fpd{L}{x^i} +
(C^\gamma_{\mu\beta}y^\mu + C^\gamma_{0\beta})\fpd{L}{y^\gamma}
-(\rho^i_0 + \rho^i_\mu y^\mu) \frac{\partial^2L}{\partial
x^i\partial y^\beta}\right),
\]
where $(g^{\alpha\beta})$ stands for the inverse matrix of
$(g_{\alpha\beta}) = \left(\frac{\partial^2 L}{\partial y^\alpha
\partial y^\beta}\right)$.

Having now seen sufficient reasons to pay particular attention
to the case of affine Lie algebroids, we come back to the construction
of Berwald-type connections associated to arbitrary $\trho$-connections.
So assume we have a $\trho$-connection on the affine Lie algebroid $\pi$
(not necessarily of pseudo-\sode\ type). As explained in Section~2,
it is then appropriate to work with the adapted basis $\{\baseH{a}, \baseV \alpha\}$
for $Sec(\pi^1)$, rather than the ``coordinate basis'' $\{\baseX{a},\baseV \alpha\}$
(here the index $\scriptsize a$ stands for either $\scriptsize 0$ or
$\scriptsize\alpha$). The following bracket relations then become useful:
\begin{equation}
\begin{array}{l}\displaystyle
[ \baseH{a}, \baseV \alpha ] =
\fpd{\Gamma_a^\delta}{y^\alpha} \baseV \delta, \\[3mm]
[\baseH{a},\baseH{b}] = C_{ab}^\delta \baseH\delta +
(C_{ab}^\delta \Gamma_\delta^\gamma +
\tilde\rho^1(\baseH{b})(\Gamma_a^\gamma)-\tilde\rho^1(\baseH{a})(\Gamma_b^\gamma))
\baseV\gamma.
\end{array} \label{hvbrackets}
\end{equation}
It will further be appropriate to write now $\subH{}$ for the projection
$\proltr E \rightarrow \pi^*\tilde E$ and likewise
define the map $\subV{}: \proltr E \rightarrow \pi^*\ov E\subset
\pi^*{\tilde E}$ by: $\subV{\mathcal Z}=\big(\trho^1(\subV
P {\mathcal Z})\big)_v$. The reason is that this will bring us in line
with notations used in \cite{Crampin,MeSa} to which the next
proposition strongly relates. Combining the horizontal and vertical
lift operations with the direct sum decomposition (\ref{decomp2})
of $Sec(\pi^*\tpi)$, it is more convenient now to think
of the following threefold decomposition of $Sec(\pi^1)$:
\begin{equation}
Sec(\pi^1)= \langle \H\cans \rangle \oplus \H{Sec(\pi^*\ov\pi)}
\oplus \V{Sec(\pi^*\ov\pi)}. \label{dec2}
\end{equation}
Note that, in the particular case of a pseudo-\sode\ connection, we have
$\H\cans=\Gamma$.

We know that any $\trho$-connection generates Berwald-type connections.
The strong point of the next result, however, is that if we assume that
$\pi$ is an affine Lie algebroid, there is a direct defining formula
for the two Berwald-type connections discussed in the preceding section.

\begin{prop} If the affine bundle $\pi$ carries an affine Lie algebroid
structure, the Berwald-type connections $\cdtilde{}$ and $\cdhattilde{}$
are determined by the following direct formulae:
\begin{eqnarray}
\cdtilde{\mathcal Z}\tilde X &=& [\subH P {\mathcal Z}, \V{\tilde
X} \subV] + [ \subV P {\mathcal Z}, \H{\tilde X} \subH] +
\trho^1(\subH P {\mathcal Z})\big(\langle
\tilde X,e^0 \rangle\big)\cans, \label{directold} \\
\cdhattilde{\mathcal Z}\tilde X &=& [ \subH P {\mathcal Z},
\V{\tilde X} \subV] + [ \subV P {\mathcal Z}, \H{\ov X} \subH] +
{\tilde\rho}^1 {\mathcal Z} \big( \langle \tilde X,e^0
\rangle\big)\cans, \label{directnew}
\end{eqnarray}
with $\ov X := \tilde X - \langle \tilde X, e^0 \rangle \cans$.
\end{prop}
\proof\, Using the properties $\subV{}\circ\subH P=0$, $\subH{}
\circ \subV P =0$, $h\circ \subH{} = \trho^1 \circ \subH{P}$, it
is easy to verify that the above expressions satisfy the
appropriate rules when the arguments are multiplied by a
function on $E$. Hence, both operators define a linear
$\trho^1$-connection on the vector bundle $\pi^*\tpi$. Next, we verify
that this connection comes from an affine $\trho^1$ connection
on $\pi^*\pi$. For that, according to a result in \cite{MeSaMa},
it is necessary and sufficient that $e^0$ (here regarded as
basic section of $\pi^*\tilde\pi$) is parallel.
We have
\begin{equation}
(\cdtilde{\mathcal Z} e^0) (\tilde X) = {\tilde\rho}^1{\mathcal Z}
(\langle\tilde X, e^0 \rangle) - \langle \cdtilde{\mathcal Z}
\tilde X , e^0\rangle, \label{cdtildee0}
\end{equation}
and similarly for $\cdhattilde{}$. In the case of $\cdtilde{}$,
we have
\[
\langle \cdtilde{\mathcal Z} \tilde X , e^0\rangle =
\langle [ \subV P {\mathcal Z}, \H{\tilde X} \subH] , e^0\rangle
+ \trho^1(\subH P {\mathcal Z})\big(\langle \tilde X,e^0 \rangle\big).
\]
Making use of the first of the bracket relations (\ref{hvbrackets}), it
is straightforward to verify that the first term on the right is
equal to $\trho^1(\subV P {\mathcal Z})\big(\langle \tilde X,e^0 \rangle\big)$,
so that the sum of both terms indeed makes the right-hand side of
(\ref{cdtildee0}) vanish. The computation for $\cdhattilde{}$ is similar.

It remains now to check that the restrictions to $Sec(\pi^*\pi)$
and $Sec(\pi^*\ov\pi)$ of (\ref{directold}) and
(\ref{directnew}) verify, respectively, the defining relations
(\ref{local1}) and (\ref{local9}) for $(\cd{},\cdov{})$ and
$(\cdhat{}, \cdhatov{})$. If we take ${\mathcal
Z}=\H{\tilde\sigma}$ and $\tilde X = \eta$, for basic $\tilde\sigma\in
Sec(\tpi)$ and $\eta\in Sec(\pi)$, then we know from
Proposition~\ref{prop5} that the bracket $[h\tilde\sigma,v\eta]$
is vertical in $TE$. As a consequence $(0,[h\tilde\sigma,v\eta])$
is vertical in $\proltr E$. But this is precisely
$[\H{\tilde\sigma},\V{\eta}]$, because the bracket of the two
projectable sections $\H{\tilde\sigma}$ and $\V{\eta}$ is by
construction (see \cite{MaMeSa}) the section
$\Big([\tilde\sigma,0],[\trho^1\H{\tilde\sigma},\trho^1\V{\eta}]\Big)$
of $\pi^1$. Therefore, $\cd{\H{\tilde\sigma}} \eta =
[\H{\tilde\sigma},\V{\eta}\subV] = \big(\trho^1(\subV
P[\H{\tilde\sigma},\V{\tilde\eta}])\big)_v
=\big(\trho^1[\H{\tilde\sigma},\V{\tilde\eta}]\big)_v
=[h\tilde\sigma,v\tilde\eta]_v$, where the Lie algebra homomorphism
provided by the anchor map $\trho^1$ has been used. Similar arguments
apply for the other operators $\cdov{}$, $\cdhat{}$ and $\cdhatov{}$ when
${\mathcal Z}$ is horizontal. In remains to look at the case ${\mathcal
Z}=\V{\ov\sigma}$ ($\ov\sigma\in Sec(\ov\pi)$). Since
$[\V{\ov\sigma}, \H{\ov\eta}]$ is vertical, it follows that
$\cdov{\V{\ov\sigma}}\ov\eta= \cdhatov{\V{\ov\sigma}}\ov\eta= 0$.
For $\tilde X = \eta$, since then $\langle\eta,e^0\rangle=1$, we
find for the first connection $\cd{\V{\ov\sigma}}\eta=0$. For the
second connection, it suffices to check (see (\ref{local7}))
that $\cdhat{\V{\ov\sigma}}\cans = 0$, and this is trivial. \qed

\section{Conclusions}

Two main objectives have been attained in this paper: we have unravelled
the mechanism by which a generalised connection over an anchored
bundle leads to a linearised connection over an appropriate prolonged
anchored bundle; we have at the same time focussed on the special features
of connections on an affine bundle, in general, and on an affine Lie
algebroid in particular. The latter subject is a completion of the
work we started in \cite{MeSaMa}. But it also ties up with the first
issue, as a generalisation of the study of Berwald-type connections
in \cite{MeSa}, where we dealt, so to speak, with the prototype
of an affine Lie algebroid, namely the first-jet extension of a
bundle fibred over $\R$, this being the geometrical arena for
time-dependent mechanics.

What are such Berwald-type connections good for?
The covariant derivative operators associated to (classical)
Berwald-type connections are those which are at the heart
of the theory of derivations of forms along the tangent
(or first-jet) bundle projection, initiated in
\cite{MaCaSaI,MaCaSaII}. These operators have proved to be
very useful tools in a number of applications concerning
qualitative features of \sode s. We mention, for example,
the characterisation of linearisability \cite{MarCar,CraMaSa}
and of separability \cite{MaCaSaIII,CaVaSa} of \sode s;
the inverse problem of Lagrangian mechanics \cite{CSMBP}; the study
of Jacobi fields and Raychaudury's equation \cite{JeriePrince}.
There is little doubt that there are similar applications ahead
for the qualitative study of pseudo-\sode s on Lie algebroids.

\end{document}